\documentclass[12pt,twoside,a4paper,leqno,onecolumn]{article}
\usepackage[english,french]{babel}
\usepackage[utf8]{inputenc}
\usepackage[OT1]{fontenc}
\usepackage{sectsty}
\usepackage{fncychap,fancybox}
\usepackage[print,panelleft,gray,paneltoc]{pdfscreen}
\usepackage{mathrsfs,amsthm,amsfonts,amssymb,amscd,amsbsy,amsmath}
\usepackage[retainorgcmds]{IEEEtrantools}
\usepackage{hyperref}
\numberwithin{equation}{section}
\newcommand {\R}{\mathbb{R}} 

\newtheorem{theo}{Th\'eor\`eme}

\newtheorem{lem}{Lemme}

\begin{document}
\pagestyle{myheadings}
\markboth{Cherrier, Hanani}{Courbure moyenne prescrite}

\thispagestyle{empty}

\centerline{\begin{Large}\bf Hypersurfaces d'un fibr\'e vectoriel Riemannien\end{Large}}

\vskip2mm

\centerline{\begin{Large}\bf \`a courbures moyennes verticale\end{Large}}

\vskip2mm

\centerline{\begin{Large}\bf et horizontale prescrites\end{Large}}

\vskip7mm

\centerline{\textbf{Pascal CHERRIER}\footnote{Adresse actuelle : Universit\'e de Paris VI, UFR 920 de Math\'ematiques, B.C. 172, 4 place Jussieu, 75252 Paris Cedex 05, France\\ \indent E-mail : cherrier@ccr.jussieu.fr}}

\vskip2mm

\centerline{et}

\vskip2mm

\centerline{\textbf{Abdellah HANANI}\footnote{Adresse actuelle : Universit\'e de Lille 1, UFR de Math\'ematiques, B\^at. M2, 59655, Villeneuve d'Ascq Cedex, France\\ \indent E-mail : abdellah.hanani@math.univ-lille1.fr}}

\vskip8mm

\hrule

\vskip4mm

\noindent\textbf{Abstract. }Let $M$ be a compact Riemannian manifold without boundary and let $E$ be a Riemannian vector bundle over $M$. If $\Sigma $ denotes the sphere subbundle of $E$, we look for embeddings of $\Sigma $ into $E$ admitting a prescribed mean curvatures of various type.

\vskip5mm

\noindent\textbf{Mots cl\'es} : Connexions, rel\`evements, courbure moyenne verticale, courbure moyenne horizontale, estimations a priori, les m\'ethodes.

\vskip3mm


\vskip4mm

\hrule

\section{Introduction}

\vskip4mm

Ce travail constitue la suite d'une \'etude portant sur la recherche d'hypersurfaces compactes d'un espace fibr\'e vectoriel Riemannien \`a courbure moyenne prescrite [4]. On d\'esigne par $(M,g)$ une vari\'et\'e Riemannienne compacte sans bord, de dimension $n\geq 1$ et $(E,{\tilde g})$ un fibr\'e vectoriel Riemannien sur $M$ de rang $m\geq 2$. On note  $\Sigma $ le fibr\'e unitaire correspondant et $E_{*}$ le fibr\'e $E$ priv\'e de la section nulle. Dans [4], on a mis en \'evidence une hypersurface de $E_{*}$ admettant une courbure moyenne \'egale \`a $K$, une fonction $\mathscr{C}^{\infty }$ strictement positive donn\'ee sur $E_{*}$, et d\'efinie comme la trace de la seconde forme fondamentale relativement \`a la m\'etrique induite par une m\'etrique Riemannienne $G$ sur $E$. La solution est donn\'ee sous la forme d'un graphe radial $\cal {Y}$ construit sur $\Sigma $, i.e. une application de $\Sigma $ dans $E_{*}$ du type $\xi \mapsto e^{u(\xi )}\xi $, o\`u $u\in \mathscr{C}^{\infty }(\Sigma )$ est une fonction inconnue qu'on prolonge \`a $E_{*}$ en la maintenant radialement constante. Les calculs ont \'et\'e \'effectu\'es dans la connexion $D$ de Sasaki.

\vskip4mm

La g\'eom\'etrie du fibr\'e ambiant permet de d\'efinir d'autres notions de courbure moyenne pour les graphes radiaux $\cal {Y}$. Une premi\`ere est d\'efinie comme suit: si $x\in M$ et si $\xi \in {\cal {Y}}_{x}=E_{x}\cap \cal {Y}$, la courbure moyenne verticale de $\cal{Y}$ au point $\xi $ est la courbure moyenne en $\xi $ de la fibre ${\cal {Y}}_{x}$ consid\'er\'ee comme hypersurface de $E_{x}$. La recherche d'un graphe radial \`a courbure moyenne verticale prescrite revient \`a la r\'esolution sur $\Sigma $ d'une \'equation elliptique d\'eg\'en\'er\'ee. Celle-ci est mise en \'evidence \`a la troisi\`eme section de cet article.

\vskip4mm

Dans le cadre euclidien, i.e. quand $M$ est r\'eduite \`a un point, les deux courbures moyenne et moyenne verticale co\"{\i}ncident. Un th\'eor\`eme de Bakelman et Kantor [2] assure en dimension $3$ l'existence d'une telle hypersurface sous la condition que la fonction $K$ d\'ecro\^{\i}t plus vite que la courbure moyenne de sph\`eres concentriques , i.e. il existe deux r\'eels $r_{1}$ et $r_{2}$ tels que $0<r_{1}\leq 1\leq r_{2}$ et
\begin{equation}
K(\xi )>\frac{m-1}{\Vert \xi \Vert }\ \mbox {si}\ \Vert \xi \Vert <r_{1},\ K(\xi
)<\frac{m-1}{\Vert \xi \Vert }\ \mbox {si}\ \Vert \xi \Vert >r_{2}
\end{equation}
jointe \`a l'hypoth\`ese de monotonicit\'e 
\begin{equation}
\frac{\partial \left[ rK(r\xi )\right]}{\partial r}\leq 0,\ \mbox{pour tout }\xi\in \Sigma.
\end{equation}
Une autre preuve, valable en toute dimension, est donn\'ee par Treibergs et Wei [9] sous les conditions pr\'ec\'edentes. L'hypoth\`ese $(1.2)$ leur a permis d'appliquer la m\'ethode de continuit\'e et leur donne l'unicit\'e \`a homoth\'etie pr\`es.

\vskip4mm

Dans le cadre des fibr\'es envisag\'es ici, le fait que l'\'equation \`a r\'esoudre soit d\'eg\'en\'er\'ee complique radicalement sa r\'esolution. D'autre part, une hypoth\`ese du type $(1.2)$ n'assure plus l'unicit\'e m\^eme \`a homoth\'etie pr\`es. Cependant, et bien qu'une r\'esolution avec une donn\'ee quelconque n'\'etait pas a priori pr\'evisible, c'est une \'etude munitieuse de cette \'equation qui am\`ene le r\'esultat suivant. 

\vskip4mm

\begin{theo}Soit $\displaystyle K\in \mathscr{C}^{\infty }(E_{*})$ une fonction partout strictement positive telle que $K(\xi )=K(\pi (\xi ))$ pour tout $\xi $, o\`u $\pi $ est la projection naturelle de E sur M. Il existe alors un graphe radial \`a courbure moyenne verticale \'egale \`a $K$.
\end{theo}

\vskip5mm

L'hypoth\`ese faite sur $K$ signifie qu'elle est le rel\`evement vertical \`a $E$ d'une fonction strictement positive de $\mathscr{C}^{\infty}(M)$. Quant \`a la preuve, elle utilise un calcul direct o\`u l'on donne explicitement une solution du probl\`eme. Remarquons que dans ce cas \'el\'ementaire, l'hypoth\`ese $(1.2)$ n'est pas satisfaite. En effet, partout dans $E_{*}$, on a : 
$$\frac{\partial \left[\rho K(\rho \xi )\right] }{\partial \rho }=K>0,\ \mbox{quel que soit }\xi \in \Sigma .$$

\vskip1mm

Par un argument de degr\'e d\'evelopp\'e dans le cadre fonctionnel $\mathscr{C}^{\infty }$ par Nagumo [8], et dont l'application repose sur l'obtention d'une estim\'ee a priori dans $\mathscr{C}^{\infty }(\Sigma )$, on d\'emontre le th\'eor\`eme suivant.

\vskip4mm

\begin{theo}Soit $\displaystyle K\in \mathscr{C}^{\infty }(E_{*})$ une fonction partout strictement positive. On suppose qu'il existe deux r\'eels $r_{1}$ et $r_{2}$, $0<r_{1}\leq 1\leq r_{2}$, tels que les in\'egalit\'es $(1.1)$ soient satisfaites. Il existe alors un graphe radial ${\cal Y}$ \`a courbure moyenne verticale donn\'ee par $K$, et tel que $r_{1}\leq \Vert \xi \Vert \leq r_{2}$ pour tout $\xi \in {\cal Y}$.
\end{theo}

\vskip5mm

A pr\'esent, on s'int\'eresse au probl\`eme de la courbure moyenne horizontale pour les graphes radiaux. Celle-ci est d\'efinie comme suit. Soit 
$$\{e_{i},e_{\alpha }\mid i=1,...,n\ {\rm et\it \ }\alpha =n+1,...,n+m\}$$
un rep\`ere mobile tangent \`a $E$, o\`u les $e_{i}$ sont des champs de vecteurs horizontaux obtenus par rel\`evement horizontal d'un rep\`ere mobile sur $M$, les 
$e_{\alpha }$ sont des champs de vecteurs verticaux et o\`u $e_{n+m}=\nu $ est le champ radial unitaire. Si $\displaystyle \xi \in E_{*}$, les $n$ vecteurs $\displaystyle e_{i}(\xi )$, $1\leq i\leq n$, forment une base du sous-espace horizontal $H_{\xi }E$ de $T_{\xi}E$. Au point ${\cal Y}(\xi )$ du graphe radial ${\cal Y}$, les $n$
vecteurs $w_{i}=D{\cal Y}(e_{i})$ forment une base du sous-espace $D{\cal
Y}(H_{\xi }E)={\cal H}_{{\cal Y}(\xi )}{\cal Y}$ de $T_{{\cal Y}(\xi)}{\cal Y}$. Notons 
${\tilde \nu }({\cal Y}(\xi ))$ l'orthogonal unitaire de ${\cal H}_{{\cal Y}(\xi )}{\cal Y}$ dans $H_{{\cal Y}(\xi )}E\oplus \R \nu ({\cal Y}(\xi ))$. Les composantes de la seconde forme fondamentale horizontale $L$ sont d\'efinies par
$$L(w_{i},w_{j})=G(D_{w_{i}}{\tilde \nu },w_{j}),\ 1\leq i,j\leq n,$$
et la courbure moyenne horizontale de ${\cal Y}$ au point ${\cal Y}(\xi )$ est alors la trace de $L$ relativement \`a la m\'etrique induite par $G$ sur ${\cal H}_{{\cal Y}(\xi )}{\cal Y}$. L'existence d'un graphe radial ${\cal Y}$ \`a courbure moyenne horizontale prescrite revient \`a r\'esoudre sur $\Sigma $ une \'equation elliptique d\'eg\'en\'er\'ee. Celle-ci est donn\'ee \`a la seconde section de cette \'etude. Elle ne peut admettre une solution si la fonction prescrite $K$ est partout strictement positive ou partout strictement n\'egative; cette remarque justifie les hypoth\`eses du
r\'esultat suivant.

\vskip4mm

\begin{theo}Soit $\displaystyle K\in \mathscr{C}^{\infty }(E_{*})$. On suppose qu'il
existe deux r\'eels $r_1$ et $r_2$, $0<r_1\leq 1\leq r_2$, tels que $\displaystyle K(\xi )>0$ si $\Vert \xi\Vert <r_1$ et $K(\xi )<0$ si $\Vert \xi \Vert >r_2$. Il existe alors un graphe radial ${\cal Y}$ \`a courbure moyenne horizontale donn\'ee par $K$ et tel que $r_1\leq \Vert \xi \Vert \leq r_2$ pour tout $\xi \in {\cal Y}$
\end{theo}

\vskip5mm

On pr\'esente la suite cette \'etude en quatre parties. La derni\`ere section est consacr\'ee \`a la preuve des th\'eor\`emes 1, 2 et 3. On y trouve aussi des exemples montrant que l'hypoth\`ese de croissance du th\'eor\`eme 2 est, dans un certain sens, la meilleure possible ainsi qu'un exemple de non unicit\'e m\^eme \`a homoth\'etie pr\`es, la condition de monotonicit\'e $(1.2)$ \'etant satisfaite. Les estimations a priori n\'ec\'essaires pour r\'esoudre dans les diff\'erents cas sont regroup\'ees \`a la quatri\`eme partie. Une mise en \'equation est pr\'esent\'ee \`a la troisi\`eme partie et, pour plus de monotonie, on donne, \`a la seconde partie de cet article, quelques rappels pr\'eliminaires et on renvoie \`a [4] pour plus de d\'etails.

\vskip6mm

\section{Rappels et notations}

\vskip4mm

\textbf{1-} Soit $(M,g)$ une vari\'et\'e Riemannienne compacte sans bord de dimension $n\geq 1$. Soient $(E,{\tilde g})$ un fibr\'e vectoriel Riemannien sur $M$ de rang $m\geq 2$, $\pi $ la projection naturelle de $E$ sur $M$ et $E_{*}$ le fibr\'e $E$ priv\'e de la section nulle. On note $\nabla $ la connexion de Levi-Civita de la vari\'et\'e $(M,g)$ et ${\tilde \nabla }$ une connexion m\'etrique sur le fibr\'e $(E,{\tilde g})$.

\vskip2mm

Soient $U$ un ouvert de $M$ muni de coordonn\'ees $(x^{i})_{1\leq i\leq n}$, $\epsilon
_{i}=\frac{\partial }{\partial x^i}$ et $\Gamma ^{k}_{ij}$ les symboles de Christoffel de $\nabla $ et $(s_{\alpha })_{n+1\leq \alpha \leq n+m}$ un rep\`ere de sections de $E$ au dessus de $U$. $\pi $ d\'esignant la projection de $E$ sur $M$, si $\xi \in \pi ^{-1}(U)$ et $x=\pi (\xi )$, on \'ecrit $\xi =y^{\alpha }s_{\alpha }(x)$ ; $(x^{i},y^{\alpha })_{i,\alpha }$ est alors un syst\`eme de coordonn\'ees sur $\pi ^{-1}(U)$. Notons $\Gamma ^{\beta }_{i\alpha }$ les symboles de Christoffel de ${\tilde \nabla }$ d\'efinis par ${\tilde \nabla}_{\epsilon _{i}}s_{\alpha }=\Gamma ^{\beta }_{i\alpha }s_{\beta }$. Le rel\`evement horizontal $e_{i}$ de $\epsilon _{i}$ est donn\'e par
\begin{equation}
e_{i}=\frac{\partial }{\partial x^i}-y^{\alpha}\Gamma ^{\beta }_{i\alpha }\frac{\partial}{\partial y^{\beta}}.
\end{equation}
Si $\displaystyle e_{\alpha }=\frac{\partial}{\partial y^{\alpha }}$, $\displaystyle {\cal S}=\{e_{i},e_{\alpha }\}_{i,\alpha }$ est un rep\`ere mobile tangent \`a $\pi ^{-1}(U)$. On d\'efinit sur $E$ une m\'etrique Riemannienne $G$ en posant
\begin{equation}
G(e_{i},e_{j})=g(\epsilon _{i},\epsilon _{j}),\ \ G(e_{\alpha },e_{\beta })={\tilde
g}(e_{\alpha},e_{\beta }),\ \ G(e_{i},e_{\alpha })=0,
\end{equation}
o\`u on identifie tout vecteur vertical \`a un point de $E$, et on consid\`ere la
connexion $D$ de Sasaki [10] d\'efinie par
\begin{equation}
D_{e_{i}}e_{j}=\Gamma ^{k}_{ij}e_{k},\ D_{e_{i}}e_{\alpha }=\Gamma ^{\beta
}_{i\alpha}e_{\beta },\ D_{e_{\alpha }}e_{i}=D_{e_{\alpha }}e_{\beta}=0.
\end{equation}
D'apr\`es [10], la connexion $D$ est compatible avec la m\'etrique $G$ et, ne
co\"{\i}ncide pas avec la connexion de Levi-Civita de $G$ ; sa torsion $T$ est non nulle et d\'epend de la courbure de ${\tilde \nabla }$. Rappellons que les composantes dans ${\cal S}$ du tenseur de courbure ${\cal R}$ de $D$ sont donn\'ees par
$$R_{dcab}=G\left((D_{e_{a}e_{b}}-D_{e_{b}e_{a}}-D_{[e_{a},e_{b}]})e_{c},e_{d}\right),\R^{d}_{cab}=G^{de}R_{ecab}$$
et un calcul direct montre, pour $1\leq i,j\leq n$ et $n+1\leq \alpha ,\beta ,\lambda ,\mu \leq n+m$, que
\begin{equation}
R^{i}_{\alpha \ \beta \ j}=R^{\lambda }_{\alpha \ \beta \ j}=R^{i}_{\alpha \ \beta \
\mu }=R^{\lambda }_{\alpha \ \beta \ \mu }\equiv 0.
\end{equation}

\vskip3mm

\noindent\textbf{2-} On note $\Sigma _{r}=\{\xi \in E\mid \Vert \xi \Vert =r\}$, $\Sigma =\{\xi \in E\mid \Vert \xi \Vert =1\}$, $\pi_{1} $ la projection naturelle du fibr\'e $\Sigma $ sur $M$, $r$ la fonction $r(\xi )=\Vert \xi \Vert $ et $\nu$ le
champ radial unitaire. Sur l'ouvert $\pi ^{-1}(U)$ muni des coordonn\'ees $(x^{i},y^{\alpha })$, $1\leq i\leq n$ et $n+1\leq \alpha \leq n+m$, le champ $\nu$ est donn\'e par
\begin{equation}
\nu =r^{-1}y^{\alpha }\frac{\partial}{\partial y^{\alpha }}.
\end{equation}
Il est normal \`a $\Sigma $ et donc l'espace tangent \`a $\Sigma $ au point $\xi \in \Sigma $ est une somme directe du sous-espace horizontal $H_{\xi }E$ de $T_{\xi
}E$ et de l'espace tangent \`a la fibre de $\Sigma $ passant par $\xi $.

\vskip2mm

Dans ce qui suit, le param\`etre $\mu _{a}$ sera \'egal \`a $0$ ou $1$ selon que
la direction $a$ est verticale ou horizontale. Fixons un rep\`ere mobile tangent \`a $E$ de la forme
$${\cal R}=\{e_{i},e_{\alpha }\mid i=1,...,n\mbox{ et }\alpha =n+1,...,n+m\},$$
o\`u les $e_{i}$ sont des champs de vecteurs horizontaux obtenus par rel\`evement horizontal d'un rep\`ere mobile $(\epsilon _{i})_{1\leq i\leq n}$ sur $M$ et o\`u les $e_{\alpha }$ sont des champs de vecteurs verticaux avec $e_{n+m}=\nu $. On notera 
$${\cal {R}}^{*}=\{\omega ^{A},\ A\leq n+m \}$$
le corep\`ere dual de ${\cal {R}}$. Appliquons $D$ \`a $e_{a}$, l'expression de $De_{a}$ dans la rep\`ere ${\cal {R}}$ nous permet d'introduire la matrice $(\omega ^{A}_{B})$ de $1$-formes d\'efinie par les \'egalit\'es
\begin{equation}
De_{A}=\omega ^{B}_{A}\otimes e_{B}.
\end{equation}
Du fait que $\nu $ est unitaire et puisque $D$ est $G$-m\'etrique, on voit que 
$$G(D_{e{_{a}}}\nu ,\nu )=0,\mbox{ pour tout }\ 1\leq a\leq n+m.$$
On en d\'eduit que
\begin{equation}
\omega ^{n+m}_{n+m}=0.
\end{equation}
D'autre part, utilisons la d\'efinition de $D$ et l'expression $(2.5)$ de $\nu $, on montre que, partout sur $\Sigma _{r}$, on a
\begin{equation}
D_{e_{a}}\nu =(1-\mu _{a})r^{-1}e_{a}\mbox{ pour }a\leq n+m-1.
\end{equation}
Reportons dans $(2.6)$, il en d\'ecoule que, partout sur $\Sigma _{r}$, on a
\begin{equation}
\omega ^{a}_{n+m}=(1-\mu _{a})r^{-1}\omega ^{a}\mbox{ pour }a\leq n+m-1
\end{equation}
et par suite, reportons $(2.7)$ et $(2.9)$ dans $(2.6)$, on obtient
\begin{equation}
D_{\nu }\nu =0.
\end{equation}
Or, $D$ \'etant $G$-m\'etrique, $\displaystyle G(D_{e_{b}}e_{a},\nu )=-G(e_{a},D_{e_{b}}\nu )$. Donc, compte tenu de $(2.8)$, partout sur $\Sigma _{r}$, on aura
\begin{equation}
\omega ^{n+m}_{a}(e_{b})=-(1-\mu _{b})r^{-1}G_{ab}\mbox{ pour }a,b\leq n+m-1.
\end{equation}
Combinons $(2.4)$, $(2.11)$ et l'\'equation de Gauss, on obtient l'expression suivante du type de composantes dans ${\cal R}$ du tenseur de courbure ${\tilde {\cal R}}$ de $\Sigma $ qui sera utilis\'e ult\'erieurement :
\begin{equation}
{\tilde R}^{j}_{\alpha \beta \gamma }={\tilde R}^{j}_{\alpha \beta i}={\tilde R}^{\gamma }_{\alpha \beta i}=0,\ n+1\leq \alpha ,\beta ,\gamma \leq n+m-1\mbox{ et }\ 1\leq i,j\leq n,
\end{equation}
et
\begin{equation}
{\tilde R}^{\lambda}_{\alpha \beta \mu }=\delta ^{\lambda }_{\beta }G_{\alpha \mu
}-\delta ^{\lambda }_{\mu }G_{\alpha \beta },\ n+1\leq \alpha ,\beta ,\lambda ,\mu
\leq n+m-1.
\end{equation} 

\vskip2mm

\noindent\textbf{3-} Soit $u\in \mathscr{C}^{2}(\Sigma )$ une fonction qu'on prolonge \`a $E_{*}$ en la maintenant radialement constante. Dans le corep\`ere ${\cal R}^*$, la
diff\'erentielle de la fonction $u$ est donn\'e par 
$$du=\sum _{a=1}^{n+m-1}D_{a}u\omega ^{a}.$$
La composante $D_{a}u$ est homog\`ene de degr\'e $(\mu _{a}-1)$. De m\^eme, on a
$$D_{ab}u=D^2u(e_{a},e_{b})=(D_{e_{a}}Du)(e_{b}).$$
D'o\`u
$$D_{ab}u=D_{e_{a}}\Big(Du(e_{b})\Big)-Du(D_{e_{a}}e_{b})$$
et on v\'erifie que la composante $\displaystyle  D_{ab}u$ est homog\`ene de degr\'e $(\mu _{a}+\mu _{b}-2)$. En particulier, on peut \'ecrire
$$D_{a\nu }u=D_{e_{a}}\left(Du(\nu )\right)-Du(D_{e_{a}}\nu )=-Du(D_{e_{a}}\nu ).$$
La derni\`ere \'egalit\'e d\'ecoule du fait que $u$ est une fonction radialement constante. Sur $\Sigma _{r}$, la relation $(2.8)$ implique alors que
\begin{equation}
D_{a\nu }u=-(1-\mu _{a})r^{-1}D_{a}u\mbox{ pour }a\leq n+m-1.
\end{equation}
Tenons compte de la relation $(2.10)$, un calcul analogue montre que
\begin{equation}
D_{\nu \nu }u=0.
\end{equation}

\vskip6mm

\section{Mise en \'equation}

\vskip4mm

On conserve les notations du dernier paragraphe et on consid\`ere l'application ${\cal Y}$ de $\Sigma $ dans $E$ telle que 
$${\cal Y}(\xi )=e^{u(\xi )}\xi ,\mbox{pour }\xi \in \Sigma ,$$
o\`u $u\in \mathscr{C}^{2}(\Sigma )$ est une fonction qu'on prolonge \`a $E_{*}$ en la
maintenant radialement constante. 

\vskip2mm

\noindent\textbf{1-} Dans ce paragraphe, on donne l'\'equation qui permet de prescrire la courbure moyenne verticale. Quand les lettres grecques sont utilis\'ees comme indice, celles ci repr\'esentent des directions verticales et d\'ecrivent l'ensemble 
$\{n+1,...,n+m-1\}$. Avec le choix ant\'erieur du rep\`ere ${\cal R}$, on voit que 
$$\displaystyle {\cal {R}}_{1}=\{e_{\alpha },\nu\mid \alpha =n+1,...,n+m-1\}$$
est un rep\`ere mobile tangent aux fibres de $E$. Notons par $\overline D$ la connexion
induite par $D$ sur les fibres de $E$. Les $(m-1)$ champs de vecteurs $e_{\alpha }$, $\alpha \in \{n+1,...,n+m-1\}$, sont tangents aux fibres de $\Sigma$, on en d\'eduit que les $(m-1)$ champs de vecteurs 
$$E_{\alpha }={\overline D}Y(e_{\alpha})=e_{\alpha }+e^{u}{\overline D}_{\alpha }u.\nu
,\ \alpha \in \{n+1,...,n+m-1\}$$
forment un rep\`ere mobile tangent aux fibres de ${\cal Y}$. Les composantes de la m\'etrique induite $h$ sont donn\'ees par 
$$h_{\alpha \beta }=G(E_{\alpha},E_{\beta})=G_{\alpha \beta }+e^{2u}{\overline D}_{\alpha }u{\overline D}_{\beta }u$$
et on v\'erifie que  
$$h^{\alpha \beta }=G^{\alpha \beta }-f^{2}e^{2u}{\overline D}^{\alpha
}u{\overline D}^{\beta }u,\ f=(1+e^{2u}{\overline D}_{\alpha }u{\overline
D}^{\alpha }u)^{-{\frac{1}{2}}}.$$
Le champ unitaire d\'efini par
$${\tilde \nu }=f(\nu -e^{u}{\overline D}^{\alpha }ue_{\alpha })$$
est normal aux fibres de ${\cal Y}$ et, compte tenu de ce choix, la courbure moyenne d'une fibre ${\cal Y}_{x}={\cal Y}\cap E_{x}$ consid\'er\'ee comme une hypersurface de $E_{x}$ est d\'efinie par
$${\cal M}_{_{{\cal Y}_{x}}}=h^{\alpha \beta }G({\overline D}_{E_{\alpha }}{\tilde \nu
},E_{\beta }).$$
Enfin, si $\xi \in {\cal Y}_{x}$, la valeur de la courbure moyenne verticale du graphe
${\cal Y}$ au point $\xi $ est 
$${\cal M}^{v}_{_{\cal Y}}(\xi )={\cal M}_{_{{\cal Y}_{x}}}(\xi ).$$
La d\'efinition de la connexion $D$ implique que $\displaystyle D_{\nu }e_{\alpha}$ est un champ de vecteurs verticaux donc, tenons compte de la relation $(2.8)$, en un point de ${\cal Y}$, on obtient
\begin{equation}
{\overline D}_{\nu }e_{\alpha }=D_{\nu }e_{\alpha}=e^{-u}e_{\alpha}.
\end{equation}
La relation $(2.14)$ se traduit alors sur chaque fibre par la suivante :
\begin{equation}
{\overline D}_{\alpha \nu }u=-e^{-u}{\overline D}_{\alpha }u\ {\rm pour\it \ }n+1\leq
\alpha \leq n+m-1.
\end{equation}
Cette derni\`ere peut \^etre \'etablise par un calcul analogue au pr\'ec\'edent. Ainsi, tenons compte de $(3.1)$, la d\'efinition de la d\'eriv\'ee covariante permet d'en d\'eduire que
\begin{equation}
{\overline D}_{\nu }({\overline D}_{\alpha }u)={\overline D}_{\nu \alpha
}u+{\overline D}u({\overline D}_{\nu }e_{\alpha})=0.
\end{equation}
D'autre part, usons du fait que $D_{\nu }u=0$ et puisque ${\overline D}_{\nu }\nu =0$, la relation $(3.3)$ permet d'\'ecrire
$${\overline D}_{E_{\alpha}}E_{\beta}={\overline D}_{e_{\alpha }}(e_{\beta }+e^{u}{\overline D}_{\beta }u.\nu )+e^{u}{\overline D}_{\alpha }u{\overline D}_{\nu }e_{\beta
}.$$
Or, la d\'efinition de la connexion $D$ montre que $D_{e_{\alpha }}e_{\beta }$ est
un champ de vecteurs verticaux. Donc, compte tenu de $(2.11)$ qui dit que $\omega ^{n+m}_{\beta }(e_{\alpha})=-e^{-u}G_{\alpha \beta }$, on obtient
$${\overline D}_{e_{\alpha }}e_{\beta }=D_{e_{\alpha }}e_{\beta}=\sum
_{n+1}^{n+m-1}\omega ^{\gamma }_{\beta }(e_{\alpha })e_{\gamma}-e^{-u}G_{\alpha \beta }\nu $$
et par suite, $(3.1)$ implique
\begin{equation}
\begin{array}{c}\displaystyle {\overline D}_{E_{\alpha }}E_{\beta }=\sum _{n+1}^{n+m-1}\omega ^{\gamma }_{\beta }(e_{\alpha })e_{\gamma }-e^{-u}G_{\alpha \beta }\nu 
+e^{u}(e_{\alpha}.{\overline D}_{\beta }u)\nu \\ \\ \displaystyle +{\overline
D}_{\beta }ue_{\alpha }+{\overline D}_{\alpha }ue_{\beta }+e^{u}{\overline
D}_{\alpha }u{\overline D}_{\beta}u\nu .\end{array}
\end{equation}
La d\'efinition de la d\'eriv\'ee covariante donne
$${\overline D}_{e_{\alpha }}({\overline D}_{\beta }u)={\overline D}_{\alpha \beta}u+\sum _{n+1}^{n+m-1}\omega ^{\gamma }_{\beta }(e_{\alpha
}){\overline D}_{\gamma }u.$$
Reportons dans $(2.4)$, la relation qui en r\'esulte peut s'\'ecrire sous la forme
$$\begin{array}{c}\displaystyle {\overline D}_{E_{\alpha }}E_{\beta }=\omega ^{\gamma}_{\beta }(e_{\alpha })E_{\gamma }+{\overline D}_{\beta }uE_{\alpha}+{\overline D}_{\alpha }uE_{\beta }\\ \\ \displaystyle +e^{-u}[-h_{\alpha \beta }+e^{2u}{\overline D}_{\alpha \beta }u]\nu \end{array}$$
et par suite, eu \'egard au fait que $G(E_{\alpha },{\tilde \nu })=0$, on obtient : 
$$G({\overline D}_{E_{\alpha }}E_{\beta},{\tilde \nu })=fe^{-u}[-h_{\alpha
\beta }+e^{2u}{\overline D}_{\alpha \beta }u].$$ 
Ainsi et puisque $G(D_{E_{\alpha }}{\tilde \nu },E_{\beta })=-G(D_{E_{\alpha }}E_{\beta
},{\tilde \nu })$, on v\'erifie qu'au point $e^{u(\xi )}\xi \in {\cal Y}_{x}$\rm , on a :
$$fe^{u}{\cal M}_{_{{\cal Y}_{x}}}=(m-1)-e^{2u}h^{\alpha \beta }{\overline
D}_{\alpha \beta }u.$$
Prenons l'image inverse sur $\Sigma _{x}$ de cette \'equation, compte tenu de l'homog\'en\'eit\'e des d\'eriv\'ees covariantes de $u$, au point $\xi \in \Sigma
_{x}$, on doit avoir :
\begin{equation}
{\overline h}^{\alpha \beta }{\overline D}_{\alpha \beta }u=(m-1)(1+{\overline
D}_{\alpha }u{\overline D}^{\alpha }u)-(1+{\overline D}_{\alpha }u{\overline D}^{\alpha
}u)^{\frac{3}{2}}e^{u}{\cal M}_{_{{\cal Y}_{x}}},
\end{equation}
o\`u ${\overline h}^{\alpha \beta}=(1+{\overline D}_{\gamma }u{\overline D}^{\gamma }u)G^{\alpha \beta }-{\overline D}^{\alpha }u{\overline D}^{\beta }u$. A pr\'esent pour toute direction verticale $\alpha $, on a : ${\overline D}_{\alpha }u=D_{\alpha }u$. D'autre part, pour tout $\alpha ,\beta \in \{ n+1,...,n+m-1\}$, utilisons la d\'efinition de la d\'eriv\'ee covariante, on peut \'ecrire
$${\overline D}_{\alpha \beta }u=e_{\alpha }({\overline D}_{\beta}u)-{\overline
D}u({\overline D}_{e_{\alpha }}e_{\beta}).$$
L'\'equation de Gauss, le fait que $u$ est une fonction radialement constante et
puisque $D_{e_{\alpha }}e_{\beta }$ est vertical impliquent que 
$${\overline D}u({\overline D}_{e_{\alpha }}e_{\beta })={\overline
D}u(D_{e_{\alpha }}e_{\beta })=Du(D_{e_{\alpha }}e_{\beta }).$$
Ainsi $\displaystyle {\overline D}_{\alpha \beta }u= D_{\alpha \beta }u$. Reportons dans l'\'equation ci-dessus, on obtient l'\'equation de la courbure moyenne verticale
\begin{equation}
\sum _{n+\leq \alpha ,\beta \leq n+m-1}B^{\alpha \beta}D_{\alpha \beta
}u=(m-1)(1+v_{1})-(1+v_{1})^{\frac{3}{2}}e^{u}{\cal M}^{v}_{_{\cal Y}}(e^u\xi ),
\end{equation}
o\`u l'on a not\'e $B^{\alpha \beta }=(1+v_{1})G^{\alpha \beta }-D^{\alpha }uD^{\beta }u$ et $v_{1}=\sum _{n+1\leq \alpha \leq n+m-1}D_{\alpha }uD^{\alpha}u$.

\vskip2mm

Remarquons que les calculs de la section pr\'ec\'edente montrent que, pour tout $r>0$, on a : $\displaystyle {\cal M}^{v}_{\Sigma _{r}}=(m-1)r^{-1}$.

\vskip5mm

\noindent\textbf{2-} Donnons ici l'\'equation de la courbure moyenne horizontale d'un
graphe radial sur $\Sigma $. Dans la suite, quand les lettres alphab\'etiques sont
utilis\'ees comme indice, celles ci repr\'esentent des directions horizontales et varient entre $1$ et $n$. Si $\displaystyle \xi \in E_{*}$, les $n$ vecteurs 
$\displaystyle e_{i}(\xi )$, $1\leq i\leq n$,, forment une base du sous-espace
horizontal $H_{\xi }E$ de $T_{\xi }E$. Au point ${\cal Y}(\xi )$ du graphe radial
$\displaystyle {\cal Y}:\xi \in \Sigma \mapsto e^{u(\xi )}\xi $, les $n$ vecteurs
$$E_{i}=D{\cal Y}(e_{i})=e_{i}+e^{u}D_{i}u\nu ,\ 1\leq i\leq n,$$
forment une base du sous-espace $D{\cal Y}(H_{\xi }E)={\cal H}_{{\cal Y}(\xi
)}{\cal Y}$ de $T_{{\cal Y}(\xi )}{\cal Y}$. Les composantes de la m\'etrique $h$ induite sur ${\cal H}_{{\cal Y}(\xi )}{\cal Y}$ sont donn\'ees par
$$h_{ij}=G(E_{i},E_{j})=G_{ij}+e^{2u}D_{i}uD_{j}u$$
et la r\'esolution de l'\'equation $h^{ij}h_{jl}=\delta ^{i}_{l}$ donne les composantes
contravariantes de celle-ci, on v\'erifie que
$$h^{ij}=G^{ij}-f^{2}e^{2u}D^{i}uD^{j}u;\ f=(1+e^{2u}D_{i}u
D^{i}u)^{-{\frac{1}{2}}}.$$
Notons ${\tilde \nu }({\cal Y}(\xi ))$ l'orthogonal unitaire de ${\cal H}_{{\cal
Y}(\xi)}{\cal Y}$ dans $H_{{\cal Y}(\xi )}E\oplus \R \nu ({\cal Y}(\xi ))$. Celui-ci est donn\'e par 
$${\tilde \nu }=f(\nu -e^{u}D^{i}ue_{i}).$$
Les composantes de la seconde forme fondamentale horizontale $L$ sont d\'efinies par
$$L(E_{i},E_{j})=G(D_{E_{i}}{\tilde \nu },E_{j}),\ 1\leq i,j\leq n,$$
et, avec le choix pr\'ec\'edent de ${\tilde \nu }$, la courbure moyenne horizontale de ${\cal Y}$ au point ${\cal Y}(\xi )$ est d\'efinie comme \'etant la trace de $L$ relativement \`a la m\'etrique induite $h$.

\vskip2mm

Tout d'abord la d\'efinition de la connexion $D$ implique que
\begin{equation}
D_{\nu }e_{i}=0.
\end{equation}
D'autre part, les relations $(2.14)$ et $(2.15)$ se traduisent par
$$D_{i\nu }u=D_{\nu \nu }u=0\mbox{ pour }1\leq i\leq n.$$
Ainsi, tenons compte de $(3.7)$, la d\'efinition de la d\'eriv\'ee covariante permet d'en d\'eduire que
\begin{equation}
D_{\nu }(D_{i}u)=D_{\nu i}u+Du(D_{\nu}e_{i})=0.
\end{equation}
Par suite, tenons compte du fait que 
$\displaystyle D_{\nu }\nu =0$ et le fait que $u$ est une constante radiale, les relations $(3.7)$ et $(3.8)$ donnent
\begin{equation}
D_{E_{i}}E_{j}=D_{e_{i}}e_{j}+(e^{u}e_{i}.D_{j}u)\nu +e^{u}D_{i}uD_{j}u\nu .
\end{equation}
Or, il d\'ecoule de la d\'efinition de $D$ que $D_{e_{i}}e_{j}$ est un champ de vecteurs horizontaux. Donc, la relation $(2.6)$ s'\'ecrit
$$D_{e_{i}}e_{j}=\sum _{1}^{n}\omega ^{k}_{j}(e_{i})e_{k}.$$
Tenons compte de cette relation, la d\'efinition de la d\'eriv\'ee covariante permet
d'\'ecrire
\begin{equation}
D_{e_{i}}(D_{j}u)=D_{ij}u+\sum _{1}^{n}\omega ^{k}_{j}(e_{i})D_{k}u.
\end{equation}
Reportons ces deux derni\`eres relations dans $(2.9)$, on obtient
$$D_{E_{i}}E_{j}=\sum _{1}^{n}\omega ^{k}_{j}(e_{i})E_{k}
+e^{u}\Big(D_{i}uD_{j}u+D_{ij}u\Big)\nu $$
et par suite, eu \'egard au fait que $G(E_{k},{\tilde \nu })=0$ pour $k\leq n$, on
obtient
$$G(D_{E_{i}}E_{j},{\tilde \nu })=fe^{u}\Big(D_{i}uD_{j}u+D_{ij}u\Big).$$
Tenons compte du fait que $G(D_{E_{i}}{\tilde \nu },E_{j})=-G(D_{E_{i}}E_{j},{\tilde \nu })$ et saturons par $h^{ij}$, on trouve l'expression suivante :
$$f^{-1}e^{u}{\cal {M}}^{h}_{_{\cal Y}}(e^u\xi )=h^{ij}\Big(-e^{2u}D_{i}uD_{j}u-e^{2u}D_{ij}u\Big).$$Prenons l'image inverse de cette \'equation sur $\Sigma $, compte
tenu de l'homogen\'eit\'e de degr\'e $0$ des d\'eriv\'ees covariantes horizontales de $u$, on obtient l'expression d\'esir\'ee :
$$f^{-1}e^{u}{\cal {M}}^{h}_{_{\cal Y}}(e^{u}\xi
)=-e^{2u}D_{i}uD^{i}u-e^{2u}h^{ij}D_{ij}u.$$
Ainsi, la recherche d'un graphe radial ${\cal Y}$ \`a courbure moyenne horizontale donn\'ee par une fonction $K$ revient \`a la r\'esolution sur $\Sigma $ de l'\'equation
elliptique d\'eg\'en\'er\'ee suivante :
\begin{equation}
\sum _{1\leq i,j\leq n}C^{ij}(u)D_{ij}u=-v_{2}-(1+e^{2u}v_{2})^{\frac{3}{2}}e^{-u}K(e^{u}\xi),
\end{equation}
o\`u on note $\displaystyle C^{ij}(u)=(1+e^{2u}v_{2})G^{ij}-e^{2u}D^{i}uD^{j}u$ et
$v_{2}=\sum _{1\leq i\leq n}D_{i}uD^{i}u$. Remarquons que le principe du maximum
implique que si $K$ est partout strictement positive ou bien partout strictement n\'egative, l'\'equation $(3.11)$ n'admet pas de solution. D'autre part, les calculs de la section pr\'ec\'edente montrent que, pour tout  $r>0$, on a : $\displaystyle {\cal M}^{h}_{\Sigma _{r}}=0$.

\vskip6mm

\section{Estimations a priori}

\vskip4mm

\begin{lem}Soient $K\in \mathscr{C}^{1}(E_{*})$ une fonction partout strictement positive et $u\in \mathscr{C}^{3}(\Sigma )$ une solution de l'\'equation
\begin{equation}
A^{ab}(u)D_{ab}u=F(\xi ,u),
\end{equation}
o\`u 
$$A^{ab}(u)=(1+v_{1}+v_{2})G^{ab}-D^{a}uD^{b}u$$
et
$$F(\xi ,u)=-(v_{1}+v_{2})^{2}+(m-1)(1+v_{1})-(1+v_{1})^{3/2}e^{u}K(e^{u}\xi )$$
telle qu'il existe deux r\'eels $r_1$ et $r_2$ v\'erifiant : $0<r_1\leq e^{u}\leq r_2$. Notons $\Sigma _{r_1,r_2}=\{\xi \in E\mid r_1\leq \Vert \xi \Vert
r_2\}$. Il existe alors une constante positive $C_{0}$ ne d\'ependant que de la g\'eom\'etrie des vari\'et\'es $(M,g)$ et $(E,{\tilde g})$, $r_1$, $r_2$, $max _{\Sigma _{r_1,r_2}}K$ et $\Vert K\Vert _{\mathscr{C}^{1}(\Sigma _{r_1,r_2})}$ telle que
$\vert Du\vert \leq C_{0}$ partout dans $\Sigma $. 
\end{lem}

\vskip5mm

\noindent\textit{D\'emonstration. }Soient $u\in \mathscr{C}^{3}(\Sigma)$ une solution
de $(4.1)$, $\ell$ un r\'eel strictement positif fix\'e ult\'erieurement et $\Gamma $ la fonctionnelle d\'efinie, sur $\Sigma $, par 
$$\Gamma (u)=(1+v)e^{\ell u},\mbox{ avec }v=v_{1}+v_{2}.$$
En un point $\xi \in \Sigma $ o\`u $\Gamma (u) $ atteint son maximum, on a
\begin{equation}
\frac{D_{a}\Gamma }{\Gamma }=\frac{D_{a}v}{1+v}+\ell D_{a}u=0
\end{equation}
et
$$A^{ab}D_{ab}\Gamma \leq 0$$
c'est-\`a-dire, tenons compte de $(4.1)$, 
\begin{equation}
A^{ab}D_{ab}v-\frac{A^{ab}D_{a}vD_{b}v}{1+v}+\ell (1+v)F\leq 0.
\end{equation}
Or,
$$A^{ab}D_{ab}v=2A^{ab}D_{a}^{c}u D_{bc}u+2A^{ab}D_{abc}uD^{c}u.$$
Une permutation de l'ordre des indices de d\'erivation covariante dans le terme des d\'eriv\'ees troisi\`emes, celle-ci g\`en\`ere des termes en torsion et en courbure, montre que
\begin{equation}
A^{ab}D_{ab}v=2A^{ab}D_{a}^{c}uD_{bc}u+2A^{ab}D_{cab}uD^{c}u+4E_{1}+E_{2},
\end{equation}
o\`u les termes $E_{1}$ et $E_{2}$ sont donn\'es par
\begin{equation}
E_{1}=-A^{ab}T^{d}_{ac}D_{bd}uD^{c}u
\end{equation}
et
\begin{equation}
E_{2}=-2A^{ab}\Big({\tilde R}^{d}\ _{bac}+D_{a}T^{d}_{bc}-T^{h}_{ac}T^{d}_{hb}\Big)
D_{d}uD^{c}u.
\end{equation}

\vskip1mm

D\'erivons une fois l'\'equation $(1.1)$ dans la direction $e_{c}$, il vient :
$$A^{ab}D_{cab}u+(D_{c}A^{ab})D_{ab}u=D_{c}F.$$
Saturons cette \'equation par $D^{c}u$ et d\'eveloppons $(D_{c}A^{ab})$, on obtient :
\begin{equation}
\begin{array}{c} \displaystyle
A^{ab}D_{cab}uD^{c}u=D_{c}FD^{c}u-D_{c}vD^{c}uD_{a}^{a}u\\ \\ \displaystyle +D_{c}\
^{a}uD^{c}uD^{b}uD_{ab}u+D_{c}\ ^{b}uD^{c}uD^{a}uD_{ab}u.\end{array}
\end{equation}
D'apr\`es $(4.1)$, on peut \'ecrire $\displaystyle D_{a}^{a}u=\frac{F+D^{a}uD^{b}uD_{ab}u}{1+v}$ et donc
$$D_{c}vD^{c}uD_{a}^{a}u=\frac{F}{1+v}D_{c}vD^{c}u+\frac{D^{a}uD^{b}uD_{ab}uD_{c}vD^{c}u}{1+v}.$$
Il en d\'ecoule que
\begin{equation}
D_{c}vD^{c}uD_{a}^{a}u=\frac{FD_{c}vD^{c}u}{1+v}+\frac{2\left(
D^{a}uD^{b}uD_{ab}u\right)^{2}}{1+v}.
\end{equation}
La d\'efinition des composantes $A^{ab}$ montre que $\displaystyle G^{ab}=\frac{A^{ab}+D^{a}uD^{b}u}{1+v}$. De ce fait, la somme des deux derniers termes du membre de droite de $(4.7)$ s'\'ecrit
$$\begin{array}{c}\displaystyle D_{c}\ ^{a}uD^{c}uD^{b}uD_{ab}u+D_{c}\ ^{b}uD^{c}uD^{a}uD_{ab}u=\frac{2\left(D^{a}uD^{b}uD_{ab}u\right)^{2}}{1+v}\\ \\ \displaystyle
+\frac{1}{1+v}A^{ad}D_{cd}uD^{c}uD^{b}u D_{ab}u+\frac{1}{1+v}A^{bd}D_{cd}uD^{c}uD^{a}u
D_{ab}u .\end{array}$$
Utilisons la relation suivante : 
\begin{equation}
D_{cd}u=D_{dc}u-T^{e}_{cd}D_{e}u,
\end{equation}
on v\'erifie que
\begin{equation}
\begin{array}{c} \displaystyle D_{c}\ ^{a}uD^{c}uD^{b}uD_{ab}u+D_{c}\
^{b}uD^{c}uD^{a}uD_{ab}u=\frac{2\left(D^{a}uD^{b}uD_{ab}u\right) ^{2}}{1+v}\\ \\
\displaystyle +\frac{1}{2(1+v)}A^{ab}D_{a}vD_{b}v+E_{3}+E_{4}. \end{array}
\end{equation}
o\`u l'on a not\'e
$$E_{4}=\frac{1}{1+v}A^{ad}T^{e}_{cd}T^{f}_{ba} D_{e}uD_{f}uD^{b}uD^{c}u$$
et
$$E_{3}=-\frac{3}{1+v}A^{ab}T^{e}_{cb}D_{e}uD^{c}uD^{d}uD_{ad}u.$$
Ce dernier peut \^etre transformer comme suit :
$$(1+v)E_{3}=-\frac{3}{2}A^{ab}T^{e}_{cb}D_{e}uD^{c}uD_{a}v,$$
ce qui, compte tenu de $(4.2)$ et la d\'efinition des composantes $A^{ab}$, permet
de v\'erifier que
\begin{equation}
E_{3}=\frac{3}{2}\ell T^{e}_{cb}D_{e}uD^{c}uD^{b}u.
\end{equation}

\vskip0mm

Reportons $(4.8)$ et $(4.10)$ dans $(4.7)$, on obtient l'\'egalit\'e suivante :
$$A^{ab}D_{cab}uD^{c}u=\frac{A^{ab}D_{a}vD_{b}v}{2(1+v)}+D_{c}FD^{c}u-\frac{FD_{a}v
D^{a}u}{1+v}+E_{3}+E_{4}.$$
Multiplions cette relation par $2$ et reportons dans $(4.4)$, l'\'egalit\'e qui en r\'esulte s'\'ecrit sous la forme : 
$$\begin{array}{c}\displaystyle A^{ab}D_{ab}v-\frac{A^{ab}D_{a}vD_{b}v}{1+v}=2D_{c}FD^{c}u-2\frac{FD_{a}vD^{a}u}{1+v}\\ \\ \displaystyle +2A^{ab}D^{c}_{a}uD_{bc}u+4E_{1}+E_{2}+2E_{3}+2E_{4}.\end{array}$$
Tenons compte de $(4.2)$, on d\'eduit de l'\'egalit\'e pr\'ec\'edente la suivante : 
\begin{equation}
\begin{array}{c}\displaystyle A^{ab}D_{ab}v-\frac{A^{ab}D_{a}vD_{b}v}{1+v}+\ell (1+v)F= 2D_{c}FD^{c}u+3\ell vF\\ \\ \displaystyle +\ell F+2A^{ab}D^{c}_{a}uD_{bc}u+4E_{1}+E_{2}+2E_{3}+2E_{4}.\end{array}
\end{equation}
D'autre part, on a : 
$$D_{c}(v_{1})D^{c}u=2(1-\mu _{a})D_{ca}uD^{a}uD^{c}u=(1-\mu
_{a})D_{a}vD^{a}u$$
de sorte que le d\'eveloppement de $D_{c}FD^{c}u$ donne 
$$\begin{array}{c}\displaystyle 2D_{c}FD^{c}u=-4vD_{c}vD^{c}u+2(m-1)(1-\mu _{a})D_{a}vD^{a}u\\ \\ \displaystyle -2(1+v_{1})^{\frac{3}{2}}e^{u}\left[v\frac{\partial [\rho K(\rho \xi )]}{\partial \rho}(e^{u}\xi)+(D_{c}K)(e^{u}\xi )D^{c}u\right]\\ \\ \displaystyle -3\sqrt{1+v_1}(1-\mu _{a})D_{a}vD^{a}ue^{u}K(e^{u}\xi).\end{array}$$
Utilisons $(4.2)$, on montre que  
\begin{equation}
\begin{array}{c}
\displaystyle 2D_{c}FD^{c}u+3\ell vF+\ell F=\ell v^{3}+3lv^{2}+(m-1)l(1+v_{2}+2v+vv_{1})\\ \\ \displaystyle -2(1+v_{1})^{\frac{3}{2}}e^{u}\left[v\frac{\partial [\rho
K(\rho \xi )]}{\partial \rho}(e^{u}\xi )+(D_{c}K)(e^{u}\xi )D^{c}u\right]\\ \\ \displaystyle -\ell (1+v+2v_2)\sqrt{1+v_1}e^{u}K(e^{u}\xi ).\end{array}
\end{equation}
Supposons que 
\begin{equation}
v(\xi )\geq C_{0}=1+8r_{2}max _{\Sigma _{r_1,r_2}}K.
\end{equation}
L'\'egalit\'e $(4.13)$ implique l'in\'egalit\'e suivante :
$$\begin{array}{c}\displaystyle 2D_{c}FD^{c}u+3\ell vF+\ell F\geq 2\left[\ell -3r_{2}\vert DK(e^{u}\xi )\vert \right]v^{2}+\ell v^{3}\\ \\ \displaystyle -6r_{2}v^{\frac{5}{2}}\left\vert\frac{\partial[\rho K(\rho \xi )]}{\partial \rho }\right\vert (e^{u}\xi)\end{array}$$
et donc, si $\ell \geq \ell _{0}=1+3r_{2}\Vert K\Vert _{\mathscr{C}^{1}(\Sigma _{r_{1},r_{2}})}$, on obtient:
$$2D_{c}FD^{c}u+3\ell vF+lF\geq \ell v^{3}-6r_{2}v^{\frac{5}{2}}\left\vert \frac{\partial [\rho K(\rho \xi )]}{\partial \rho }\right\vert (e^{u}\xi ).$$
Reportons cette in\'egalit\'e dans $(4.12)$, compte tenu de $(4.14)$ et puisque $\vert E_{2}\vert \leq C_{1}v^{2}$, o\`u $C_{1}$ est fonction de $\Vert {\cal T}\Vert _{\infty }$, $\Vert D{\cal T}\Vert _{\infty }$ et $\Vert {\tilde {\cal R}}\Vert _{\infty }$, $\vert E_{4}\vert \leq C'_{1}v^{2}$ et, d'apr\`es $(4.11)$,
$\displaystyle \vert E_{3}\vert \leq C_{2}\ell v^{\frac{3}{2}}$, o\`u $C'_{1}$ et $C_{2}$ ne d\'ependent que de $\Vert {\cal T}\Vert _{\infty }$, on v\'erifie que 
$$\begin{array}{c}\displaystyle A^{ab}D_{ab}v-\frac{A^{ab}D_{a}vD_{b}v}{1+v}+\ell (1+v)F\geq -C_{3}\ell v^{2}+\ell v^{3}+2A^{ab}D^{c}_{a}uD_{bc}u
\\ \\ \displaystyle +4E_{1}-6r_{2}v^{\frac{5}{2}}\left\vert \frac{\partial
[\rho K(\rho \xi )]}{\partial \rho }\right\vert (e^{u}\xi ),\end{array}$$
o\`u \it $C_{3}$ \rm est une constante positive ne d\'ependant que des constantes $C_{1}$, $C'_{1}$ et $C_{2}$. Ainsi, compte tenu de $(4.3)$, on aboutit \`a l'in\'egalit\'e suivante: 
$$\left[\ell (v-C_{3})-6r_{2}v^{\frac{1}{2}}\left\vert
\frac{\partial [\rho K(\rho \xi)]}{\partial \rho }\right\vert
(e^{u}\xi )\right]v^{2}+2A^{ab}D^{c}_{a}uD_{bc}u+4E_{1}\leq 0.$$
De sorte que si
\begin{equation}
v(X)\geq C'_{3}=2C_{3}+\left[12r_{2}\left\Vert \frac{\partial [\rho K(\rho
\xi )]}{\partial \rho }\right\Vert _{\mathscr{C}^{0}(\Sigma _{r_1,r_2})}\right]^{2},
\end{equation}
on obtient :
\begin{equation}
A^{ab}D^{c}_{a}uD_{bc}u+2E_{1}\leq 0.
\end{equation}

\vskip1mm

Maintenant, on d\'eveloppe le carr\'e suivant :
$$\begin{array}{c}\displaystyle K=G^{ab}G^{cd}\left(D_{ad}u-\varepsilon D_{a}uD^{e}uD_{ed}u-T^{l}_{ae}D^{e}uG_{ld}\right)\\ \\ \displaystyle \left(D_{bc}u-\epsilon D_{b}uD^{e}uD_{ec}u-T^{l}_{be}D^{e}uG_{lc}\right)\end{array}$$
avec $\displaystyle \varepsilon =\frac{1}{1+v}$ pour voir que
$$\begin{array}{c}\displaystyle (1+v)K=A^{ab}D^{c}_{a}uD_{bc}u+ (1+v)G^{ab}G_{cd}T^{c}_{ae}D^{e}uT^{d}_{bf}D^{f}u\\ \\ \displaystyle +2E_{1}-(1+v)^{-1}D^{a}uD^{b}uD_{ac}uD_{b}^{c}u.\end{array}$$
Il existe, alors, une constante positive $C_{4}$ fonction de $\Vert {\cal T}\Vert _{\infty }$, telle que
\begin{equation}
A^{ab}D^{c}_{a}uD_{bc}u+2E_{1}\geq -C_{4}v(1+v)+(1+v)^{-1}D^{a}uD^{b}uD_{ac}uD_{b}^{\
c}u.
\end{equation}
Utilisons $(4.9)$, on montre que
$$\begin{array}{c}\displaystyle D^{a}uD^{b}uD_{ac}uD_{b}^{\ c}u=D^{a}uD^{b}uD_{ca}uD^{c}_{\ b}u-2D^{a}uD^{b}uD^{c}_{\ a}uT^{e}_{bc}D_{e}u\\ \\
\displaystyle +G^{cd}D^{a}uD^{b}uT^{e}_{ac}D_{e}uT^{f}_{bd}D_{f}u\end{array}$$
et par suite
$$D^{a}uD^{b}uD_{ac}uD_{b}^{\ c}u={\frac{1}{4}}D^{a}vD_{a}v-(D^{b}uD^{c}v)T^{e}_{bc}D_{e}u+G^{cd}D^{a}uD^{b}uT^{e}_{ac}D_{e}uT^{f}_{bd}D_{f}u.$$
Ainsi, tenons compte de $(4.2)$ et $(4.14)$, il existe une constante positive $C_{5}$ fonction de $\Vert {\cal T}\Vert _{\infty }$, telle que
$$D^{a}uD^{b}uD_{ac}uD_{b}^{\ c}u\geq \frac{\ell ^2}{4}(1+v)^{2}v-C_{5}\ell (1+v)^{\frac{3}{2}}v.$$
Reportons cette in\'egalit\'e dans $(4.17)$, on en d\'eduit que
$$A^{ab}D^{c}_{a}uD_{bc}u+2E_{1}\geq (4^{-1}\ell ^{2}-C_{4})v(1+v)-\ell C_{5}v\sqrt{1+v}.$$ 
Une in\'egalit\'e qui, compte tenu de $(4.16)$, implique
$$(4^{-1}\ell ^{2}-C_{4})\sqrt{1+v}-\ell C_{5}\leq 0.$$
De sorte que pour $\ell $ assez grand, il existe une constante positive $C_{6}$ telle que $v(\xi )\leq C_{6}$. Ainsi, compte tenu de $(4.14)$ et $(4.15)$, on voit que
$$v(\xi )\leq C_{7}=max(C_{0},C'_{3},C_{6}).$$
La d\'efinition de la fonctionnelle $\Gamma $ montre que partout dans $\Sigma $, on a: 
$$v\leq \left(1+C_{7}\right)\left(\frac{r_2}{r_1}\right)^{\ell }.$$
Le lemme est prouv\'e. 

\vskip5mm

\begin{lem} Conservons les notations du lemme pr\'ec\'edent. Soient $K\in \mathscr{C}^{1}(E_{*})$ et $u\in \mathscr{C}^{3}(\Sigma )$ une solution de l'\'equation suivante :
\begin{equation}
A^{ab}D_{ab}u=-(v_{1}+v_{2})^{2}-v_{2}-(1+e^{2u}v_{2})^{\frac{3}{2}}e^{-u}K(e^{u}\xi),
\end{equation}
o\`u l'on a not\'e
$$A^{ab}(u)=(1+e^{2u}v_{1}+e^{2u}v_{2})G^{ab}-e^{2u}D^{a}uD^{b}u.$$
On suppose qu'il existe deux r\'eels $r_1$ et $r_2$ tels que $0<r_1\leq e^{u}\leq r_2$. Alors il existe une constante positive $C_0$ ne d\'ependant que de la g\'eom\'etrie des vari\'et\'es $(M,g)$ et $(E,{\tilde g})$, $r_1$, $r_2$
et $\Vert K\Vert _{\mathscr{C}^{1}(\Sigma _{r_1,r_2})}$ telle que $\vert Du\vert \leq C_{0}$ partout dans $\Sigma $. 
\end{lem}

\vskip5mm

\noindent\textit{D\'emonstration. }Soient $u\in \mathscr{C}^{3}(\Sigma )$ une solution de $(4.18)$, $\ell $ un r\'eel strictement positif fix\'e ult\'erieurement et $\Gamma $ la fonctionnelle d\'efinie par 
$$\Gamma (u)=ve^{\ell u}.$$
En un point $\xi \in \Sigma $ o\`u $\Gamma (u) $ atteint son maximum, supposons que $v(\xi )\geq 1$. Signalons qu'il suffit de majorer $v(\xi )$, la d\'efinition de $\Gamma $ permet de conclure. Les calculs qui suivent seront \'evalu\'es au point $\xi $ et l'on a
\begin{equation}
\frac{D_{a}\Gamma }{\Gamma }=\frac{D_{a}v}{v}+\ell D_{a}u=0
\end{equation}
et
$$A^{ab}D_{ab}\Gamma \leq 0.$$
Tenons compte de l'\'equation $(4.18)$ satisfaite par $u$ et notons
$$F=-v^{2}-v_{2}-(1+e^{2u}v_{2})^{\frac{3}{2}}e^{-u}K(e^{u}\xi ),$$
la derni\`ere in\'egalit\'e s'\'ecrit comme suit :
$$A^{ab}D_{ab}v-\frac{A^{ab}D_{a}vD_{b}v}{v}+\ell vF\leq 0.$$
Compte tenu de $(4.19)$ et puisque $A^{ab}D_{a}uD_{b}u=v$, on aboutit \`a
\begin{equation}
A^{ab}D_{ab}v-\ell ^{2}v^{2}+\ell vF\leq 0.
\end{equation}
Or
$$A^{ab}D_{ab}v=2A^{ab}D_{a}^{\ c}uD_{bc}u+2A^{ab}D_{abc}uD^{c}u.$$
Une permutation de l'ordre des indices de d\'erivation covariante dans le terme en d\'eriv\'ees troisi\`emes, celle-ci g\`en\`ere des termes en torsion et en courbure, montre que 
\begin{equation}
A^{ab}D_{ab}v=2A^{ab}D_{a}^{\ c}uD_{bc}u+2A^{ab}D_{cab}uD^{c}u+4E_{1}+E_{2},
\end{equation}
o\`u les termes $E_{1}$ et $E_{2}$ sont donn\'es par
\begin{equation}
E_{1}=-A^{ab}T^{d}_{ac}D_{bd}uD^{c}u
\end{equation}
et 
$$E_{2}=-2A^{ab}\Big({\tilde R}^{d}\ _{bac}+D_{a}T^{d}_{bc}-T^{h}_{ac}T^{d}_{hb}\Big) D_{d}uD^{c}u.$$
En particulier, il existe une constante positive $C_{1}$ ne d\'ependant que de $\Vert u\Vert _{\infty }$, $\Vert {\cal T}\Vert _{\infty }$, $\Vert D{\cal T}\Vert _{\infty }$ et $\Vert {\tilde {\cal R}}\Vert _{\infty }$ telle que
\begin{equation}
\vert E_{2}\vert \leq C_{1}v^{2}.
\end{equation}

\vskip0mm

D\'erivons une fois l'\'equation $(4.18)$ dans la direction $e_{c}$ et saturons l'\'equation ainsi obtenue par $D^{c}u$, il vient
\begin{equation}
A^{ab}D_{cab}uD^{c}u+(D_{c}A^{ab})D^{c}uD_{ab}u=D_{c}FD^{c}u.
\end{equation}
Le d\'eveloppement de $(D_{c}A^{ab})$ montre que
$$\begin{array}{c}\displaystyle \left(D_{c}A^{ab}\right)D^{c}uD_{ab}u=e^{2u}D_{c}vD^{c}uG^{ab}D_{ab}u+2e^{2u}v^{2}G^{ab}D_{ab}u\\ \\ \displaystyle -2e^{2u}vD^{a}uD^{b}uD_{ab}u-e^{2u}\left[D_{c}\ ^{a}uD^{c}uD^{b}uD_{ab}u+D^{a}uD^{c}uD_{c}\ ^{b}uD_{ab}u\right].\end{array}$$
Une \'egalit\'e qui, compte tenu de l'\'equation $(4.18)$ satisfaite par $u$ et d'o\`u l'on extrait la valeur de $G^{ab}D_{ab}u$, s'\'ecrit sous la forme :
$$\begin{array}{c}\displaystyle \left(D_{c}A^{ab}\right)D^{c}uD_{ab}u=\frac{e^{2u}F}{1+e^{2u}v}\left(D_{a}vD^{a}u+2v^{2}\right)\\ \\ \displaystyle -e^{2u}\left[D_{c}\ ^{a}uD^{c}uD^{b}uD_{ab}u+D^{a}uD^{c}uD_{c}\ ^{b}uD_{ab}u\right] \\ \\ \displaystyle +\frac{e^{2u}}{1+e^{2u}v}\left(e^{2u}D_{c}vD^{c}u-2v\right)D^{a}uD^{b}uD_{ab}u .\end{array}$$
Reportons dans $(4.24)$, on obtient la relation suivante :
$$\begin{array}{c}\displaystyle A^{ab}D_{cab}uD^{c}u= D_{a}FD^{a}u-\frac{e^{2u}F}{1+e^{2u}v}\left(D_{a}vD^{a}u+2v^{2}\right)\\ \\ \displaystyle +e^{2u}\left[D_{c}\ ^{a}uD^{c}uD^{b}uD_{ab}u+ D^{a}uD^{c}uD_{c}\ ^{b}uD_{ab}u\right]\\ \\ \displaystyle
-\frac{e^{2u}}{1+e^{2u}v}\left(e^{2u}D_{c}vD^{c}u-2v\right)D^{a}uD^{b}uD_{ab}u.\end{array}$$
Or, $\displaystyle D_{a}v=2D^{b}uD_{ab}u$. Donc, la pr\'ec\'edente s'\'ecrit sous la forme :
\begin{equation}
\begin{array}{c}\displaystyle A^{ab}D_{cab}uD^{c}u=D_{a}FD^{a}u-\frac{e^{2u}F}{1+e^{2u}v}\left(D_{a}vD^{a}u+2v^{2}\right)\\ \\ \displaystyle +e^{2u}\left[2^{-1}D_{c}\ ^{a}uD^{c}uD_{a}v+D^{a}uD^{c}uD_{c}\ ^{b}uD_{ab}u\right]\\ \\ \displaystyle -\frac{e^{4u}}{2(1+e^{2u}v)}D_{a}vD^{a}uD_{b}vD^{b}u+\frac{ve^{2u}}{1+e^{2u}v}D_{a}vD^{a}u.\end{array}
\end{equation}
Usons de la commutation
\begin{equation}
D_{cd}u=D_{dc}u-T^{e}_{cd}D_{e}u,
\end{equation}
on peut \'ecrire l'\'egalit\'e $(4.25)$ sous la forme :
\begin{equation}
\begin{array}{c}\displaystyle A^{ab}D_{cab}uD^{c}u =D_{a}FD^{a}u-\frac{e^{2u}F}{1+e^{2u}v}\left(D_{a}vD^{a}u+2v^{2}\right)+2^{-1}e^{2u}D^{a}vD_{a}v\\ \\
\displaystyle -\frac{e^{4u}}{2(1+e^{2u}v)}D^{a}uD^{b}uD_{a}vD_{b}v
+\frac{ve^{2u}}{1+e^{2u}v}D_{a}vD^{a}u+2^{-1}E_{3}+E_{4},\end{array}
\end{equation}
o\`u les termes $E_{3}$ et $E_{4}$ sont donn\'es par 
$$E_{3}=-3e^{2u}T^{e}_{ab}D_{e}uD^{a}uD^{b}v$$
et
$$E_{4}=e^{2u}G^{ad}D^{c}uD^{d}uT^{e}_{ac}T^{f}_{bd}D_{e}uD_{f}u.$$
Remarquons qu'il existe une constante positive $C_{2}$ ne d\'ependant que de 
$\Vert u\Vert _{\infty }$ et $\Vert {\cal T}\Vert _{\infty }$ telle que
\begin{equation}
\vert E_{4}\vert \leq C_{2}v^{2}.
\end{equation}
D'une part, tenons compte de $(4.19)$, on voit que le terme $E_{3}$ s'\'ecrit sous la forme :
$$E_{3}=3\ell e^{2u}vT^{e}_{ab}D_{e}uD^{a}uD^{b}u.$$
D'o\`u, l'existence d'une constante positive $C_3$ ne d\'ependant que de $r_2$ et $\Vert {\cal T}\Vert _{\infty }$ telle que
\begin{equation}
\vert E_{3}\vert \leq C_{3}\ell v^{\frac{5}{2}}.
\end{equation}
D'autre part, l'\'egalit\'e $(4.27)$ donne la suivante :
$$\begin{array}{c}\displaystyle A^{ab}D_{cab}uD^{c}u=D_{a}FD^{a}u+(\ell -2)\frac{e^{2u}v^{2}F}{1+e^{2u}v}\\ \\ \displaystyle +\frac{\ell (\ell -2)e^{2u}v^{3}}{2(1+e^{2u}v)}+2^{-1}E_{3}+E_{4}.\end{array}$$
Multiplions cette relation par $2$ et reportons dans $(4.21)$, on obtient
l'\'egalit\'e suivante :
$$\begin{array}{c}\displaystyle A^{ab}D_{ab}v=2D_{a}FD^{a}u+2(\ell -2)\frac{e^{2u}v^{2}F}{1+e^{2u}v}+2A^{ab}D_{a}^{\ c}uD_{bc}u \\ \\ \displaystyle
+\frac{\ell (\ell -2)e^{2u}v^{3}}{1+e^{2u}v}+4E_{1}+E_{2}+E_{3}+2E_{4}.\end{array}$$
Par suite, l'in\'egalit\'e $(4.20)$ se transforme comme suit :
\begin{equation}
\begin{array}{c}\displaystyle -\frac{\ell v^{2}(\ell +2e^{2u}v)}{1+e^{2u}v}+\frac{(3\ell -4)e^{2u}v+\ell }{1+e^{2u}v}vF\\ \\ \displaystyle +2D_{a}FD^{a}u+2A^{ab}D_{a}^{\ c}uD_{bc}u+4E_1+E_2+E_3+2E_4\leq 0.\end{array}
\end{equation}
A pr\'esent, remarquons que
$$D_{a}v_{2}D^{a}u=2\mu _{b}D^{a}uD^{b}uD_{ab}u=\mu _{a}D^{a}uD_{a}v$$
de sorte que le d\'eveloppement de $D_{a}FD^{a}u$ donne 
$$\begin{array}{c}\displaystyle D_{a}FD^{a}u=-2vD_{a}vD^{a}u
-\mu _{a}D_{a}vD^{a}u+(1+e^{2u}v_{2})^{\frac{3}{2}}ve^{-u}K(e^{u}\xi )\\ \\ \displaystyle -(1+e^{2u}v_{2})^{\frac{3}{2}}\left[v\frac{\partial [K(\rho \xi )]}{\partial \rho }(e^{u}X)+e^{-u}(D_{a}K)(e^{u}\xi )D^{a}u\right]\\ \\ \displaystyle -\frac{3}{2}(1+e^{2u}v_{2})^{\frac{1}{2}}(\mu _{a}D^{a}vD_{a}u+2vv_{2})e^{u}K(e^{u}\xi ),\end{array}$$
ce qui, compte tenu de $(4.19)$, montre que
$$\begin{array}{c}\displaystyle D_{a}FD^{a}u=2\ell v^{3}+\ell vv_{2}+(1+e^{2u}v_{2})^{\frac{3}{2}}ve^{-u}K(e^{u}\xi)\\ \\ \displaystyle -(1+e^{2u}v_{2})^{\frac{3}{2}}\left[v\frac{\partial [K(\rho \xi )]}{\partial \rho}(e^{u}\xi)+e^{-u}(D_{a}K)(e^{u}\xi )D^{a}u\right]\\ \\ \displaystyle +\frac{3}{2}(\ell -2)(1+e^{2u}v_{2})^{\frac{1}{2}}vv_{2}e^{u}K(e^{u}\xi ).\end{array}$$
Reportons dans $(4.30)$, tenons compte de l'expression de $F$, l'in\'egalit\'e qui en r\'esulte s'\'ecrit :
$$\begin{array}{c}\displaystyle (\ell +4)v^{3}-2\ell v^{2}+(4-3\ell )vv_{2}-\frac{\ell (\ell -2)v^{2}}{1+e^{2u}v}+\frac{2\ell (\ell -2)v(v^{2}-v_{2})}{1+e^{2u}v}\\ \\ \displaystyle -2(1+e^{2u}v_{2})^{\frac{3}{2}}\left[v\frac{\partial [K(\rho \xi )]}{\partial \rho}(e^{u}\xi )+e^{-u}(D_{a}K)(e^{u}\xi )D^{a}u\right]\\ \\ \displaystyle -(\ell -2)(1+e^{2u}v_{2})^{\frac{1}{2}}e^{-u}K(e^{u}\xi )\left[v+\frac{2e^{2u}vv_{1}}{1+e^{2u}v}\right]\\ \\ \displaystyle +2A^{ab}D_{a}^{\ c}uD_{bc}u+4E_1+E_2+E_3+2E_4\leq 0.\end{array}$$
Or $v_{2}\leq v$ et $1\leq v$. Donc, supposons $\ell \geq 2$, on obtient : 
$$\begin{array}{c}\displaystyle \ell v^{3}-\ell (\ell +5)v^{2}+2A^{ab}D_{a}^{\ c}uD_{bc}u+4E_1+E_2+E_3+2E_4\\ \\ \displaystyle -2(1+e^{2u}v_{2})^{\frac{3}{2}}\left[v\frac{\partial [K(\rho \xi )]}{\partial \rho }(e^{u}\xi )+e^{-u}(D_{a}K)(e^{u}\xi)D^{a}u\right]\\ \\ \displaystyle -(\ell -2)(1+e^{2u}v_{2})^{\frac{1}{2}}e^{-u}K(e^{u}\xi)\left[v+\frac{2e^{2u}vv_{1}}{1+e^{2u}v}\right]\leq 0.\end{array}$$
Ainsi, compte tenu de $(4.23)$, $(4.28)$ et $(4.29)$, il existe une constante positive $C_{4}$, fonction de $C_{1}$, $C_{2}$, $C_{3}$, $\Vert K\Vert _{\infty }$, $\Vert DK\Vert _{\infty }$ et $\Vert u\Vert _{\infty }$, telle que
\begin{equation}
\ell v^{3}-C_{4}\ell ^{2}v^{\frac{5}{2}}+2A^{ab}D_{a}^{\ c}uD_{bc}u+4E_{1}\leq 0.
\end{equation}
Pour conclure \`a partir de $(4.31)$, on d\'eveloppe le carr\'e suivant :
$$\begin{array}{c}\displaystyle K=G^{ab}G^{cd}\left(D_{ad}u-\varepsilon
D_{a}uD^{e}uD_{ed}u-T^{l}_{ae}D^{e}uG_{ld}\right)\\ \\ \displaystyle \left(D_{bc}u-\varepsilon D_{b}uD^{e}uD_{ec}u T^{l}_{be}D^{e}uG_{lc}\right)\end{array}$$
avec $\displaystyle \varepsilon =\frac{e^{2u}}{1+e^{2u}v}$ pour voir d'abord que 
$$\begin{array}{c}\displaystyle (1+e^{2u}v)K=A^{ab}D^{c}_{a}uD_{bc}u (1+e^{2u}v)G^{ab}G_{cd}T^{c}_{ae}D^{e}uT^{d}_{bf}D^{f}u\\ \\ \displaystyle +2E_{1}-\frac{e^{2u}}{1+e^{2u}v}D^{a}uD^{b}uD_{ac}uD_{b}^{c}u\end{array}$$
et de la positivit\'e de $K$ et par un raisonnement analogue \`a celui de la fin de la preuve du lemme 1, on montre qu'il existe deux constantes positives $C_{5}$ et $C_{6}$ ne d\'ependant que de $\Vert u\Vert _{\infty }$ et $\Vert {\cal T}\Vert _{\infty }$ telles que
$$A^{ab}D^{c}_{a}uD_{bc}u+2E_{1}\geq (4^{-1}\ell ^{2}-C_{5})v^{2}-C_{6}\ell v^{\frac{3}{2}}.$$
Reportons dans $(4.31)$ et simplifions par $\displaystyle v^{\frac{3}{2}}$, on obtient :
\begin{equation}
\ell v(\sqrt{v}-\ell C_{4})+(2^{-1}\ell ^{2}-2C_{5})\sqrt{v}\leq 2\ell C_{6}.
\end{equation}
Posons $\displaystyle \ell =2\sqrt{C_{5}+1}\geq 2$, on voit que $\displaystyle
2^{-1}\ell ^{2}-2C_{5}=1$ et si 
$$\displaystyle v(\xi )\geq (\ell C_{4}+1)^{2},$$
alors, d'apr\`es $(4.32)$, $v(\xi )\leq (2\ell C_{6})^{2}$. En conclusion, $v(\xi)\leq max\left((\ell C_{4}+1)^{2},(2\ell C_{6})^{2}\right)$ et le lemme est prouv\'e.

\vskip6mm

\section{Preuve des r\'esultats}
\subsection{Preuve du th\'eor\`eme 1}

\vskip4mm

\noindent On consid\`ere la fonction $\displaystyle u=-\log (K)$. L'hypoth\`ese faite sur $K$ implique que 
$$u=u\circ \pi \ \mbox{ partout dans }E_{*}.$$
Notons par $\pi _{*}$ l'application lin\'eaire tangente de $\pi $ de sorte que 
$$Du(Y)=Du\left( \pi _{*}(Y)\right) \mbox{ pour tout }Y\in TE.$$
Or, en tout point $\xi \in E_{*}$, l'espace vertical de $T_{\xi }E$ n'est autre que le noyau de $\pi _{*\mid T_{\xi }E}$. Il en d\'ecoule que, pour toute direction verticale
$\alpha $,
$$D_{\alpha }u=0.$$
Et comme $D_{e_{\alpha }}e_{\beta }$ est un champ de vecteurs vertical. La d\'efinition de la d\'eriv\'ee covariante montre que, pour $\alpha ,\
\beta \in \{n+1,...,n+m-1\}$, on a:
$$D_{\alpha \beta }u=0.$$
Reportons dans l'\'equation $(3.6)$, on voit que la courbure moyenne verticale du graphe $\xi \in \Sigma \mapsto e^{u(\xi )}\xi $ est donn\'ee par
$${\cal M}^{v}(e^{u}\xi )=K(\xi )=K(e^{u}\xi ),\ \mbox{ pour tout }\xi \in \Sigma .$$
D'o\`u le r\'esultat.

\vskip6mm

\subsection{Preuve du th\'eor\`eme 2}

\vskip4mm

Rappelons que $K\in \mathscr{C}^{\infty }(E_{*})$ est une fonction partout strictement positive donn\'ee et qu'on suppose l'existence de deux r\'eels $r_1$ et $r_2$ dans $]0,+\infty [$, $r_1\leq 1\leq r_2$, v\'erifiant
\begin{equation}
\left\{\begin{array}{ll}
\displaystyle K(\xi )>\frac{m-1}{\Vert \xi \Vert }&\mbox{si } \Vert \xi \Vert <r_{1}\\ \\ \displaystyle K(\xi )<\frac{m-1}{\Vert \xi \Vert}&\mbox{si }\Vert \xi \Vert >r_{2}.
\end{array}\right.
\end{equation}
Chercher un graphe radial ${\cal Y}$ sur $\Sigma $ \`a courbure moyenne verticale donn\'ee par $K$ revient, d'apr\`es les calculs de la seconde section, \`a r\'esoudre dans $\Sigma $ l'\'equation suivante : 
\begin{equation}
\sum _{n+1\leq \alpha ,\beta \leq n+m-1}B^{\alpha \beta }(u)D_{\alpha
\beta }u=(m-1)(1+v_{1})-(1+v_{1})^{\frac{3}{2}}e^{u}K(e^{u}\xi).
\end{equation}

\vskip2mm

\noindent La d\'emonstration qu'on donne ici est une adaptation de celle du th\'eor\`eme 3 dans [4]. On note $\displaystyle \Sigma '=\{\xi \in E\mid r_1\leq \Vert \xi \Vert \leq r_2\}$, $\nu$ le champ radial unitaire et $r$ la fonction $r(\xi )=\Vert \xi \Vert $. Soit
$\displaystyle w\in \mathscr{C}^{\infty }(\Sigma ')$, on d\'esigne par ${\tilde w}$ la restriction de $w$ \`a $\Sigma $ que l'on prolonge en une constante radiale et on note 
$$\begin{array}{l}\displaystyle v_{1}(w)=(1-\mu _{a})D_{a}wD^{a}w,\
v_{2}(w)=\mu _{a}D_{a}wD^{a}w,\ v(w)=v_{1}(w)+v_{2}(w),\\ \\
\displaystyle A^{ab}(w)=\left(1+r^{2}v_{1}({\tilde w})+r^{(1-\mu _{a})(1-\mu _{b})}v_{2}({\tilde w})\right)G^{ab}-r^{2-\mu _{a}-\mu _{b}}D^{a}{\tilde
w}D^{b}{\tilde w},\\ \\ \displaystyle F(w)=(m-1)[1+r^{2}v_{1}({\tilde
w})]-\left[1+r^{2}v_{1}({\tilde w})\right]^{3/2}e^{{\tilde w}}
K\left(e^{{\tilde w}}\frac{\xi }{\Vert \xi \Vert }\right), \\ \\
B(w)=(m-1)\left[1+r^{2}v_{1}{\tilde w})+v_{2}({\tilde w})\right]-r^{2}v_{1}({\tilde w}).\end{array}$$
Pour $t\in [0,1]$ et $w\in \mathscr{C}^{\infty }(\Sigma ')$, d\'esignons par $T_{t}(w)=u_{t}$ l'unique solution du probl\`eme de Neumann 
\begin{equation}
\left\{\begin{array}{ll}\displaystyle D_{\nu \nu }u+A^{ab}(w)r^{2-\mu
_{a}\mu _{b}}D_{ab}u-u=&t\left[-{\tilde w}+F(w)\right]+D_{\nu \nu }w\\ &+rB(w)D_{\nu }w\mbox{ dans } \Sigma '\\ \displaystyle D_{\nu }u=0\ \mbox{ sur }\partial \Sigma '.&\end{array}\right.
\end{equation}
La preuve de l'existence d'une solution $ u_{t}\in C^{2,\alpha }(\Sigma ')$ de $(5.3)$ est tr\`es classique, on renvoie \`a [1] et [5], voir aussi [6], quant \`a l'unicit\'e, elle d\'ecoule du principe du maximum. L'\'equation $(5.3)$ \'etant elliptique, on en d\'eduit qu'il existe une constante positive $C$ telle que
$$\displaystyle \Vert u_{t}\Vert _{\mathscr{C}^{2,\alpha }(\Sigma
')}\leq C.$$
Or toutes les donn\'ees sont de classe $\mathscr{C}^{\infty}$, par r\'egularit\'e, la solution $u_{t}$ l'est aussi. En particulier, \'etant donn\'ee une partie ${\cal B}$ born\'ee de $\mathscr{C}^{\infty }(\Sigma ')$, il existe ${\tilde u_{t}}\in \mathscr{C}^{\infty }(\Sigma ')$ ainsi
qu'une suite de r\'eels positifs $(C_{k})_{k\geq 0}$ telle que, quel que soit $(t,w)\in [0,1]\times {\cal B}$ et pour tout entier $k\geq 0$, on ait : 
\begin{equation}
\Vert u_{t}\Vert _{\mathscr{C}^{k}(\Sigma ')}\leq C_{k}.
\end{equation}
Il en d\'ecoule que l'op\'erateur $T_{t}$ d\'efini sur $\mathscr{C}^{\infty }(\Sigma ')$ est compact. Il en est de m\^eme de l'op\'erateur $T$ d\'efini de $[0,1]\times \mathscr{C}^{\infty }(\Sigma ')$ vers $\mathscr{C}^{\infty }(\Sigma ')$ par 
\begin{equation}
T(t,w)=w-T_{t}w
\end{equation}
et l'on veut montrer que l'\'equation
\begin{equation}
T(t,u)=0,\ {\rm pour\ \it }t\in [0,1],
\end{equation}
admet une solution dans $\mathscr{C}^{\infty }(\Sigma ')$. Si une telle solution existe, celle-ci, not\'ee $u_{t}$, est une constante radiale. En effet, $u_{t}$ v\'erifie le syst\`eme suivant :
\begin{equation}
\left\{ \begin{array}{c}\displaystyle A^{ab}(u_{t})r^{2-\mu
_{a}\mu _{b}}D_{ab}u_{t}-u_{t}=t\left[-{\tilde u_{t}}+F(u_{t})\right] +rB(u_{t})D_{\nu }u_{t}\mbox{ dans }\Sigma '\\ \\ \displaystyle D_{\nu }u_{t}=0\mbox{ sur }
\partial \Sigma ',\end{array}\right.
\end{equation}
o\`u l'on a not\'e par ${\tilde u_{t}}$ le prolongement en une constante radiale de la restriction \`a $\Sigma $ de $u_{t}$.

\vskip2mm

D\'erivons radialement l'\'equation $(5.7)$ et multiplions par $r$ l'\'equation ainsi obtenue. Du fait que $D_{\nu }r=1$ et puisque ${\tilde u_{t}}$ est une constante radiale, il en d\'ecoule que
\begin{equation}
\begin{array}{c}\displaystyle rD_{\nu}\left[A^{ab}(u_{t})\right]r^{2-\mu _{a}\mu
_{b}}D_{ab}u_{t}+(2-\mu _{a}\mu _{b})r^{2-\mu _{a}\mu _{b}}A^{ab}(u_{t})D_{ab}u_{t}\\ \\ \displaystyle +A^{ab}(u_{t})r^{2-\mu _{a}\mu _{b}}rD_{\nu ab}u_{t}-rD_{\nu
}u_{t}=trD_{\nu}F(u_{t})\\ \\ \displaystyle +B(u_{t})(rD_{\nu }u_{t}+r^{2}D_{\nu \nu }u_{t})+r^{2}D_{\nu }u_{t}D_{\nu}B(u_{t}).\end{array}
\end{equation}
D'abord, d'apr\`es les calculs de la premi\`ere section,
\begin{equation}
D_{e_{a}}(r\nu )=(1-\mu _{a})e_{a},
\end{equation}
o\`u $\mu _{a}$ vaut $1$ ou $0$ selon que la direction $e_{a}$ est horizontale ou verticale. Ainsi, usons de la d\'efinition de la d\'eriv\'ee covariante, on peut\'ecrire:
$$D_{e_{a}}(rD_{\nu }u)=D^{2}u(e_{a},r\nu)+Du\left(D_{e_{a}}(r\nu )\right)=rD_{\nu }(D_{a}u)+(1-\mu _{a})D_{a}u.$$
On en d\'eduit que
$$rD_{\nu}(D_{a}u)=D_{e_{a}}(rD_{\nu }u)-(1-\mu _{a})D_{a}u$$
et, en particulier, pour une fonction $u$ radialement constante, on a
\begin{equation}
\begin{array}{c}\displaystyle rD_{\nu }v_{1}(u)=2(1-\mu _{a})D_{a}(rD_{\nu }u)
D^{a}u-2v_{1}(u)=-2v_{1}(u),\\ \\ \displaystyle rD_{\nu }v_{2}(u)=2\mu
_{a}{\tilde D}_{a}(rD_{\nu }u) {\tilde D}^{a}u+2v_{2}(u)rD_{\nu }u=0.\end{array}
\end{equation}
Combinons $(5.10)$ avec le fait que ${\tilde u_{t}}$ est une fonction radialement constante, on obtient:
$$rD_{\nu }\Big[r^{2}v_{1}({\tilde u_{t}})+v_{2}({\tilde
u_{t}})\Big]=0.$$
Ainsi, l'\'equation $(5.8)$ se r\'eduit \`a la suivante :
\begin{equation}
\begin{array}{c}\displaystyle (2-\mu _{a}\mu _{b})r^{2-\mu _{a}\mu
_{b}}A^{ab}(u_{t})D_{ab}u_{t}+A^{ab}(u_{t})r^{2-\mu _{a}\mu _{b}}rD_{2n,ab}u_{t}\\ \\ \displaystyle +(1-\mu _{a})(1-\mu _{b})r^{3-\mu _{a}-\mu _{b}}v_{2}({\tilde u_{t}})G^{ab}D_{ab}u_{t}-rD_{\nu }u_{t}\\ \\ \displaystyle =-tr\left[r^{2}v_{1}({\tilde
u_{t}})+v_{2}({\tilde u_{t}})\right]^{2}+B(u_{t})(rD_{\nu }u_{t}+r^{2}D_{\nu \nu }u_{t}).\end{array}
\end{equation}
D'autre part, la d\'efinition de la d\'eriv\'ee covariante nous permet d'\'ecrire que
$$D_{ab}(rD_{\nu }u)=D^{2}(rD_{\nu }u)(e_{a},e_{b})=D_{e_{a}}\left[D(rD_{\nu
}u)(e_{b})\right]-D(rD_{\nu }u)(D_{e_{a}}e_{b}).$$ 
D'o\`u 
$$\begin{array}{c}\displaystyle D_{ab}(rD_{\nu}u)=e_{a}\left[D_{e_{b}}\left(Du(r\nu
)\right)\right]-\left(D_{e_{a}}e_{b}\right)Du(r\nu )\\ \\ \displaystyle
=D_{e_{a}}(D^{2}u(e_{b},r\nu )+D_{e_{a}}Du\left(D_{e_{b}}(r\nu )\right)\\ \\ \displaystyle -D^{2}u\left(D_{e_{a}}e_{b},r\nu \right)-Du\left(D_{D_{e_{a}}e_{b}}(r\nu )\right)\end{array}$$
et par suite 
$$\begin{array}{c}\displaystyle D_{ab}(rD_{\nu }u)=D^{3}u(e_{a},e_{b},r\nu )+D^{2}u\left(e_{b},D_{e_{a}}(r\nu )\right)+D^{2}u\left(e_{a},D_{e_{b}}(r\nu )\right)\\ \\
\displaystyle +Du\left(D_{e_{a}}(D_{e_{b}}(r\nu ))\right)-Du\left(D_{D_{e_{a}}e_{b}}(r\nu )\right).\end{array}$$
Or, la relation $(5.10)$ permet d'\'ecrire que
$$Du\left(D_{e_{a}}(D_{e_{b}}(r\nu ))\right)=(1-\mu _{b})Du\left(D_{e_{a}}e_{b}\right)$$
et 
$$Du\left(D_{D_{e_{a}}e_{b}}(r\nu )\right)=(1-\mu _{b})Du\left(D_{e_{a}}e_{b}\right).$$
Ainsi, usons \`a nouveau de $(5.10)$, on obtient : 
$$D_{ab}(rD_{\nu }u)=D^{3}u(e_{a},e_{b},r\nu
)+(2-\mu _{a}-\mu _{b})D^{2}u(e_{a},e_{b})$$
et d'apr\`es la d\'efinition de la connexion $D$ et celle de la d\'eriv\'ee covariante,
on voit que
\begin{equation}
D_{ab}(rD_{\nu }u)=rD_{\nu }(D_{ab}u)+(2-\mu _{a}-\mu _{b})D_{ab}u.
\end{equation}
Eu \'egard \`a l'expression $(2.4)$ donnant l'expression de la courbure de $D$ et dont d\'ecoule l'\'egalit\'e $\displaystyle D_{ab\nu }u=D_{\nu ab}u$. Reportons $(5.12)$ dans $(5.11)$, on obtient:
\begin{equation}
\begin{array}{c}
\displaystyle A^{ab}(u_{t})r^{2-\mu _{a}\mu _{b}}D_{ab}(rD_{\nu }u_{t})+(\mu _{a}+\mu _{b}-\mu _{a}\mu _{b})r^{2-\mu _{a}\mu _{b}}A^{ab}(u_{t})D_{ab}u_{t}\\ \\ \displaystyle +(1-\mu _{a})(1-\mu _{b})r^{3-\mu _{a}-\mu _{b}}v_{2}({\tilde u_{t}})G^{ab}D_{ab}u_{t}-rD_{\nu }u_{t}\\ \\ \displaystyle =B(u_{t})(rD_{\nu }u_{t}+r^{2}D_{\nu \nu
}u_{t})-tr\left[r^{2}v_{1}({\tilde u_{t}})+v_{2}({\tilde u_{t}})\right]^{2}.\end{array}
\end{equation}

\vskip0mm

Rappelons que l'\'equation de Gauss et les calculs de la seconde section impliquent que pour deux fonctions $u_{1},u_{2}\in \mathscr{C}^{2}(\Sigma ')$ ayant les m\^emes valeurs sur $\Sigma $, on a:
\begin{equation}
D_{ab}u_{1}=D_{ab}u_{2}+(1-\mu _{a})G_{ab}D_{\nu }(u_{1}-u_{2}),\; a,b\leq n+m-1.
\end{equation}
Or ${\tilde u_{t}}=(u_{t})_{\vert \Sigma }$. Donc, tenons compte de $(5.14)$ et du fait que $D_{\nu }u_{t}$ est nul sur $\Sigma $, l'\'equation $(5.13)$ montre que, partout dans $\Sigma $, on a :
\begin{equation}
(\mu _{a}+\mu _{b}-\mu _{a}\mu _{b})A^{ab}(u_{t})D_{ab}{\tilde
u_{t}}+(1-\mu _{a})v_{2}({\tilde u_{t}})D^{a}_{a}{\tilde
u_{t}}=-t\left[v_{1}({\tilde u_{t}})+v_{2}({\tilde u_{t}})\right]^{2}.
\end{equation} 
La fonction ${\tilde u_{t}}$ \'etant une constante radiale, il en d\'ecoule que celle-ci v\'erifie l'\'equation
\begin{equation}
\begin{array}{c}\displaystyle (\mu _{a}+\mu _{b}-\mu _{a}\mu _{b})r^{2-\mu _{a}-\mu
_{b}}A^{ab}(u_{t})D_{ab}{\tilde u_{t}}+(1-\mu _{a})r^{2(1-\mu _{a})}v_{2}({\tilde
u_{t}})D^{a}_{a}{\tilde u_{t}}\\ \\ \displaystyle =-t\left[r^{2}v_{1}({\tilde
u_{t}})+v_{2}({\tilde u_{t}})\right]^{2}
\end{array}
\end{equation}
partout dans $\Sigma '$. D'autre part, restreignons $(5.7)$ \`a $\Sigma $, l'usage de $(5.14)$ montre que la fonction ${\tilde u_{t}}$ v\'erifie, partout dans $\Sigma $, l'\'equation suivante :
$$A^{ab}(u_{t})D_{ab}{\tilde u_{t}}-{\tilde u_{t}}=t\left[-{\tilde u_{t}}+F({\tilde u_{t}})\right]$$
et par soustraction de l'\'equation $(5.15)$ de cette derni\`ere, on voit que, partout dans $\Sigma$, on a: 
$$\begin{array}{c}\displaystyle \sum _{n+1\leq \alpha , \beta \leq n+m-1}
B^{\alpha \beta }(u_{t})D_{\alpha \beta }{\tilde u_{t}}-{\tilde u_{t}}=t\left[-{\tilde
u_{t}} +(m-1)\left(1+v_{1}({\tilde u_{t}})\right)\right] \\ \\
\displaystyle -t\Big[1+v_{1}({\tilde u_{t}})\Big]^{\frac{3}{2}}e^{{\tilde
u_{t}}}K(e^{{\tilde u_{t}}}\xi ).\end{array}$$
La fonction ${\tilde u_{t}}$ \'etant une constante radiale, il en d\'ecoule que celle-ci v\'erifie l'\'equation
\begin{equation}
\begin{array}{c}\displaystyle (1-\mu _{a})(1-\mu _{b})r^{2-\mu _{a}-\mu _{b}}A^{ab}(u_{t})D_{ab}{\tilde u_{t}}-(1-\mu _{a})r^{3(1-\mu _{a})}v_{2}({\tilde u_{t}})D^{a}_{a}{\tilde u_{t}}\\ \\ \displaystyle -{\tilde u_{t}}=t\left\{-{\tilde u_{t}}+(m-1)\left[1+v_{1}({\tilde u_{t}})\right]-\left[1+v_{1}({\tilde u_{t}})\right]^{3/2}e^{{\tilde u_{t}}}K(e^{{\tilde u_{t}}}\xi )\right\}
\end{array}
\end{equation}
partout dans $\Sigma '$. Multiplions $(5.16)$ par $r$ et sommons l'\'equation ainsi obtenue avec $(5.17)$, montre que ${\tilde u_{t}}$ est une autre solution de $(5.7)$. En particulier, on v\'erifie que 
$$\begin{array}{c}\displaystyle A^{ab}({\tilde u_{t}})r^{2-\mu _{a}\mu _{b}}D_{ab}({\tilde u_{t}}-u_{t})-({\tilde u_{t}}-u_{t})=rB({\tilde u_{t}})D_{\nu }({\tilde
u_{t}}-u_{t})\ \mbox{ dans }\Sigma '\\ \\ \displaystyle D_{\nu }({\tilde u_{t}}-u_{t})=0\ \mbox{ sur }\partial \Sigma '.\end{array}$$
Le principe du maximum implique que ${\tilde u_{t}}=u_{t}$ partout dans $\Sigma
_{1,2}$. Ceci montre que $u_{t}$ est une constante radiale solution dans $\Sigma $ du syst\`eme suivant:
\begin{equation}
\begin{array}{c} \displaystyle A^{ab}(u_{t})D_{ab}u_{t}-u_{t}=t\left[-u_{t} +F(u_{t})\right]\\ \\ \displaystyle \sum _{n+1\leq \alpha , \beta \leq n+m-1} B^{\alpha \beta }(u_{t})D_{\alpha \beta}u_{t}-u_{t}=t\left[-u_{t}+(m-1)\left(1+v_{1}(u_{t})\right)\right] \\ \\ \displaystyle -t\left[1+v_{1}(u_{t})\right]^{3/2}e^{u_{t}}K(e^{u_{t}}\xi ).\end{array}
\end{equation}
A pr\'esent, on montre que $u_{t}$ est estim\'ee \`a priori dans $\mathscr{C}^{0}(\Sigma )$. Une majoration a priori se d\'eduit imm\'ediatement du principe du maximum. Soit $\xi \in \Sigma $ un point o\`u $u_{t}$ atteint son maximum. Si $u_{t}(\xi )>\log (r_2)$, l'hypoth\`ese de croissance $(5.1)$ faite sur $K$ combin\'ee avec $(5.18)$ implique qu'au point $\xi $, on aura 
$$0\geq D_{a}^{a}u_{t}=(1-t)u_{t}+t(m-1)-te^{u_{t}}K(e^{u_{t}}\xi )>0$$
ce qui constitue une contradiction, eu \'egard au fait que $\log (r_2)\geq 0$.
La minoration $u_{t}\geq \log (r_1)$ s'obtient par analogie en consid\'erant un point o\`u $u_{t}$ atteint son minimum.  

\vskip2mm

La fonction $u_{t}$ v\'erifiant la premi\`ere \'equation dans $(5.18)$, d'apr\`es le lemme 1 elle est estim\'ee a priori dans  $\mathscr{C}^{1}(\Sigma )$. La th\'eorie classique des \'equations uniform\'ement elliptiques, cf. [5] ou [7], assure l'existence d'un r\'eel $c_{1}>0$ tel que
$$\Vert u_{t}\Vert _{\mathscr{C}^{1,\alpha }(\Sigma )}<c_{1}$$
D'autre part, les in\'egalit\'es de Schauder implique l'existence d'une constante positive $C$ telle que
$$\Vert u_{t}\Vert _{\mathscr{C}^{2,\alpha }(\Sigma )}<C\left[\Vert u_{t}\Vert
_{\mathscr{C}^{0}(\Sigma )}+\Vert F\Vert _{\mathscr{C}^{0,\alpha }(\Sigma
')}\right]$$
et donc il existe une constante positive $c_{2}$ telle que
$$\Vert u_{t}\Vert _{\mathscr{C}^{2,\alpha }(\Sigma )}<c_{2}.$$
Supposons donc que que pour tout $s\leq k$, pour un $k\geq 2$, on ait:
$$\Vert u_{t}\Vert _{\mathscr{C}^{s,\alpha }(\Sigma )}<c_{s}.$$
L'\'equation obtenue en d\'erivant $(k-1)$ fois la premi\`ere \'equation dans $(5.18)$ s'\'ecrit localement sous la forme
$$A^{ab}D_{ab}(D_{i_{1}i_{2}...i_{k-1}}u_{t})=H_{i_{1}i_{2}...i_{k-1}},$$
o\`u le second membre $\displaystyle H_{i_{1}i_{2}...i_{k-1}}$ ne d\'epend que des
d\'eriv\'ees covariantes de $u_{t}$ d'ordre $\leq k$. Il est donc born\'e dans 
$\mathscr{C}^{0,\alpha }(\Sigma ')$, et il en d\'ecoule, d'apr\`es les in\'egalit\'es de Schauder, que $\displaystyle \Vert D^{(k-1)}u_{t}\Vert _{\mathscr{C}^{2,\alpha }(\Sigma )}<C_{k+1}$, et par suite
$$\Vert u_{t}\Vert _{\mathscr{C}^{k+1,\alpha }(\Sigma )}<c_{k+1}.$$
On vient ainsi d'\'etablir par r\'ecurrence l'existence de r\'eels positifs $a_{k}$ tels que
\begin{equation}
\Vert u_{t}\Vert _{\mathscr{C}^{k}(\Sigma ')}<a_{k}.
\end{equation}
Eu \'egard au fait que $u_{t}$ est radialement constante. On note $B$ la pseudo-boule d\'efinie par
$$B=\{ u\in \mathscr{C}^{\infty
}(\Sigma ')\mid \Vert u\Vert _{\mathscr{C}^{k}(\Sigma ')}<a_{k}\}.$$
Compte tenu de $(5.19)$, l'\'equation $(5.6)$ n'admet pas de solution $u$ sur le bord $\partial B$ de $B$. La d\'eformation $T$ est donc une homotopie compacte sur le bord de $B_{R}$. En cons\'equence, d'apr\`es le th\'eor\`eme de Nagumo [8], le d\'egr\'e de $T$ en $0$ relativement \`a $B$ ne d\'epend pas de $t$. Ainsi, pout tout $t\in [0,1]$, on a : 
\begin{equation}
d(T(t,.),0,B)=d(T(0,.),0,B)=\gamma .
\end{equation}
Or, pour $t=0$, la fonction $u_{0}=0$ est l'unique solution de $(5.6)$ et l'on montre ais\'ement que pour $\displaystyle w\in \mathscr{C}^{\infty }(\Sigma ')$,  $(d_{u_{0}}T_{0})(w)=u$, o\`u $u$ est l'unique solution du probl\`eme suivant :
\begin{equation}
\left\{\begin{array}{cl}
\displaystyle D_{\nu \nu }u+r^{2-2\mu _{a}}D^{a}_{a}u+\mu
_{a}\log (r)D^{a}_{a}u-u= D_{\nu \nu }w+(m-1)rD_{\nu }w&\mbox{dans }\Sigma
'\\ \\ \displaystyle D_{\nu }u=0&\mbox{sur }\partial \Sigma '.\end{array}\right.
\end{equation}
A pr\'esent, si $w\in \ker[\mbox{id}-(d_{u_{0}}T_{0})]$, celle-ci sera solution du probl\`eme suivant :
$$\left\{\begin{array}{cl}\displaystyle r^{2-2\mu _{a}}D^{a}_{a}w+\mu _{a}\log (r)D^{a}_{a}w-w=(m-1)rD_{\nu }w&\mbox{dans }\ \Sigma '\\ \\ \displaystyle D_{\nu }w=0&\mbox{sur }\partial \Sigma '.\end{array}\right.$$ 
Un raisonnement analogue \`a celui qui pr\'ec\`ede montre que $w$ est une constante radiale. Celle-ci \'etant identiquement nulle sur $\Sigma $, donc partout nulle et par suite
\begin{equation}
\ker[\mbox{id}-(d_{u_{0}}T_{0})]=\{0\}.
\end{equation}
Or, d'apr\`es sa d\'efinition $(5.21)$ et en raisonnant comme auparavant, on v\'erifie que $d_{u_{0}}T_{0}$ est un op\'erateur compact, et donc, tenons compte de $(5.22)$, le th\'eor\`eme d'alternative de Fredholm implique que l'op\'erateur $\displaystyle \mbox{id}-(d_{u_{0}}T_{0})$ est inversible. En cons\'equence $0$ est un point r\'egulier pour $\mbox{id}-T_{0}$. Ainsi, dans $(5.20)$, $\gamma =\pm 1$ et en particulier, l'op\'erateur $T_{1}$ admet un point fixe qui est une constante radiale de classe $\mathscr{C}^{\infty}(\Sigma )$ solution de $(5.2)$.

\vskip6mm

\subsection{Preuve du th\'eor\`eme 3}

\vskip4mm

\noindent Soit $K\in \mathscr{C}^{\infty }(E_{*})$ une fonction donn\'ee telle qu'il existe deux r\'eels strictement positifs $r_1$ et $r_2$, $r_1\leq 1\leq r_2$,
v\'erifiant
\begin{equation}
K(\xi )>0\mbox{ si }\Vert \xi \Vert <r_1\mbox{ et }K(\xi )<0\mbox{ si }\Vert \xi
\Vert >r_2.
\end{equation}
La recherche d'un graphe radial sur $\Sigma $ \`a courbure  moyenne horizontale donn\'ee par $K$, revient, d'apr\`es les calculs du second paragraphe de la seconde section, \`a r\'esoudre sur $\Sigma $ l'\'equation suivante :
\begin{equation}
\sum _{1\leq i,j\leq n}C^{ij}(u)D_{ij}u=-v_{2}-(1+e^{2u}v_{2})^{\frac{3}{2}}e^{-u}K(e^{u}\xi).
\end{equation}

\vskip1mm

Conservons les notations de la preuve du th\'eor\`eme 2 mais, avec la convention qu'ici
$$\begin{array}{l}\displaystyle A^{ab}(w)=\left(1+r^{2+\mu _{a}\mu _{b}}e^{2{\tilde w}}v_{1}({\tilde w})+e^{2{\tilde w}}v_{2}({\tilde w})\right)G^{ab}-r^{2-\mu _{a}-\mu
_{b}}e^{2{\tilde w}}D^{a}{\tilde w}D^{b}{\tilde w},\\ \\ \displaystyle F(w)=-v_{2}({\tilde w})-[1+e^{2{\tilde w}}v_{2}({\tilde w})]^{\frac{3}{2}}e^{-{\tilde
w}}K\left(e^{{\tilde w}}\frac{\xi }{\Vert \xi \Vert }\right)-r[r^{2}v_{1}({\tilde w})+v_{2}({\tilde w})]^{2}\\ \\ \displaystyle B(w)=(m-1)[1+r^{2}e^{2{\tilde w}}v_{1}{\tilde w})+e^{2{\tilde w}}v_{2}({\tilde w})]-r^{2}e^{2{\tilde w}}v_{1}({\tilde w}).\end{array}$$
Pour $t\in [0,1]$ et $w\in \mathscr{C}^{\infty }(\Sigma ')$, d\'esignons par $T_{t}(w)=u_{t}$ l'unique solution du probl\`eme de Neumann
\begin{equation}
\left\{ \begin{array}{cl} \displaystyle D_{\nu \nu }u+A^{ab}(w)r^{\delta _{ab}}D_{ab}u-u=&t\Big[-{\tilde w}+F(w)\Big]+D_{\nu \nu }w\\ &+rB(w)D_{\nu }w\mbox{ dans }\Sigma '\\ \displaystyle D_{\nu }u=0\mbox{ sur }\partial \Sigma ',& \end{array}\right.
\end{equation}
o\`u $\delta _{ab}=3-\mu _{a}-\mu _{b}-\mu _{a}\mu _{b}$. L'\'equation $(5.25)$ \'etant elliptique et toutes les donn\'ees sont de classe $\mathscr{C}^{\infty }$, un raisonnement analogue \`a celui qui pr\'ec\`ede implique l'existence de r\'eels positifs $(C_{k})_{k\geq 0}$ telle que, quel que soit $t\in [0,1]$ et pour tout entier $k\geq 0$, on ait : 
$$\Vert {\tilde u_{t}}\Vert _{\mathscr{C}^{k}(\Sigma ')}\leq C_{k}.$$
Il en d\'ecoule que l'op\'erateur $T_{t}$ est compact, il en est de m\^eme de l'op\'erateur $T$ d\'efini de $[0,1]\times \mathscr{C}^{\infty }(\Sigma ')$ vers $\mathscr{C}^{\infty }(\Sigma ')$ par
$$T(t,w)=w-T_{t}w$$
et l'on veut montrer que l'\'equation 
\begin{equation}
T(t,u)=0,\mbox{ pour }t\in [0,1],
\end{equation}
admet une solution dans $\mathscr{C}^{\infty }(\Sigma ')$. Si une telle solution existe, celle-ci, not\'ee $u_{t}$, est une constante radiale. En effet, $u_{t}$ v\'erifie le syst\`eme suivant :
\begin{equation}
\left\{\begin{array}{cl} \displaystyle A^{ab}(u_{t})r^{\delta _{ab}}D_{ab}u_{t}-u_{t}=t\left[-{\tilde u_{t}}+F(u_{t})\right]+rB(u_{t})D_{\nu }u_{t}&\mbox{dans }\Sigma '\\ \\ \displaystyle D_{\nu }u_{t}=0&\mbox{sur }\partial \Sigma ',\end{array}\right.
\end{equation}
o\`u l'on a not\'e par ${\tilde u_{t}}$ le prolongement en une constante radiale de
la restriction \`a $\Sigma $ de $u_{t}$.

\vskip2mm

D\'erivons radialement $(5.27)$ et multiplions par $r$ l'\'equation ainsi obtenue. Les calculs men\'es pour prouver le th\'eor\`eme 2 montrent 

\begin{equation}
\begin{array}{c} \displaystyle A^{ab}(u_{t})r^{\delta _{ab}}D_{ab}(rD_{\nu }u_{t})+\delta _{ab}r^{\delta _{ab}}A^{ab}(u_{t})D_{ab}u_{t}\\ \\ \displaystyle +\mu _{a}\mu _{b}r^{5-\mu _{a}-\mu _{b}}e^{2{\tilde u_{t}}}v_{1}({\tilde u_{t}})G^{ab}D_{ab}u_{t}-rD_{\nu }u_{t}\\ \\ \displaystyle =B(u_{t})(rD_{\nu }u_{t}+r^{2}D_{\nu \nu }u_{t})-tr
\Big[r^{2}v_{1}({\tilde u_{t}})+v_{2}({\tilde u_{t}})\Big]^{2}. \end{array}
\end{equation}

\vskip0mm

Or ${\tilde u_{t}}=(u_{t})_{\vert \Sigma }$. Donc, tenons compte de $(5.14)$, preuve du th\'eor\`eme 2, et du fait que $D_{\nu }u_{t}$ est nul sur $\Sigma $, l'\'equation $(5.28)$ montre que, partout dans $\Sigma $, on a :
\begin{equation}
(1-\mu _{a}\mu _{b})A^{ab}(u_{t})D_{ab}{\tilde u_{t}}+\mu _{a}e^{2{\tilde
u_{t}}}v_{1}({\tilde u_{t}})D^{a}_{a}{\tilde u_{t}}=-t\left[v_{1}({\tilde u_{t}})+v_{2}({\tilde u_{t}})\right]^{2}.
\end{equation}
D'autre part, restreignons $(5.27)$ \`a $\Sigma $, l'usage de $(5.14)$, preuve du th\'eor\`eme 2, montre que la fonction ${\tilde u_{t}}$ v\'erifie, partout dans $\Sigma $, l'\'equation suivante :
$$A^{ab}(u_{t})D_{ab}{\tilde u_{t}}-{\tilde u_{t}}=t\left[-{\tilde u_{t}}+F({\tilde
u_{t}})\right]$$
et par soustraction de l'\'equation $(5.29)$ de cette derni\`ere, on voit que
\begin{equation}
\begin{array}{c} \displaystyle \sum _{1\leq i,j\leq n}C^{ij}(u_{t})D_{ij}u_{t}-u_{t}=t\left[-u_{t}-v_{2}(u_{t})\right] \\ \\ \displaystyle -t\left[1+e^{2u_{t}}v_{2}(u_{t})\right]^{\frac{3}{2}}e^{-u_{t}}K(e^{u_{t}}\xi ). \end{array}
\end{equation}
La fonction ${\tilde u_{t}}$ \'etant une constante radiale, une combinaison des \'equations $(5.29)$ et $(5.30)$, montre que ${\tilde u_{t}}$ est une autre solution de $(5.27)$. En particulier, on v\'erifie que
$$\left\{ \begin{array}{cl}\displaystyle A^{ab}({\tilde u_{t}})r^{\delta
_{ab}}D_{ab}({\tilde u_{t}}-u_{t})-({\tilde u_{t}}-u_{t})=rB({\tilde u_{t}})D_{\nu
}({\tilde u_{t}}-u_{t})&\mbox{dans }\Sigma '\\ \\ \displaystyle D_{\nu }({\tilde u_{t}}-u_{t})=0&\mbox{sur }\partial \Sigma '.\end{array}\right.$$
Le principe du maximum montre que ${\tilde u_{t}}=u_{t}$ partout dans $\Sigma '$. Ceci montre que $u_{t}$ est une fonction radialement constante solution dans $\Sigma$ du syst\`eme :
\begin{equation}
\left\{ \begin{array}{l} \displaystyle A^{ab}(u_{t})D_{ab}u_{t}-u_{t}=t\left[-u_{t} +F(u_{t})\right]\\ \\\displaystyle C^{ij}(u_{t})D_{ij}u_{t}-u_{t}=t\left[-u_{t}-v_{2}(u_{t})\right]-t\left[1+e^{2u_{t}}v_{2}(u_{t})\right]^{\frac{3}{2}}e^{-u_{t}}K(e^{u_{t}}\xi ).\end{array}\right.
\end{equation}
Une telle solution est estim\'ee \`a priori dans $\mathscr{C}^{0}(\Sigma )$. En effet, une majoration a priori se d\'eduit imm\'ediatement du principe du maximum. Soit $\xi \in \Sigma $ un point o\`u $u_{t}$ atteint son maximum. Si $u_{t}(\xi )>\log (r_2)$, l'hypoth\`ese de croissance $(5.23)$ faite sur $K$ combin\'ee avec $(5.31)$ implique qu'au point $\xi $, on aura  
$$0\geq D_{a}^{a}u_{t}=(1-t)u_{t}-te^{-u_{t}}K(e^{u_{t}}\xi )>0$$
ce qui constitue une contradiction, eu \'egard au fait que $\log (r_{2})\geq 0$. La minoration $u_{t}\geq \log (r_{1})$ s'obtient par analogie en consid\'erant un point o\`u $u_{t}$ atteint son minimum.

\vskip2mm

La fonction $u_{t}$ v\'erifiant la premi\`ere \'equation dans $(5.31)$, le lemme 2 implique, via la th\'eorie classique des \'equations uniform\'ement elliptiques, que 
$$\Vert u_{t}\Vert _{\mathscr{C}^{1,\alpha }(\Sigma )}<Cste.$$
Ainsi en raisonnant par r\'ecurrence comme auparavant, compte tenu du fait
que $u_{t}$ est radialement constante, on conclut \`a l'existence de r\'eels positifs $a_{k}$ tels que 
$$\Vert u_{t}\Vert _{\mathscr{C}^{k}(\Sigma ')}<a_{k}.$$
Ceci nous permet d'appliquer le m\^eme argument topologique que celui utilis\'e pour prouver le th\'eor\`eme 2 pour affirmer l'existence d'un point fixe de $T_{1}$ et l'on voit, d'apr\`es $(5.31)$, que celui-ci est solution de $(5.24)$.

\vskip6mm

\subsection{Remarques}

\vskip4mm

\noindent\textbf{1-} Expliquons ici en quel sens l'hypoth\`ese de croissance du th\'eor\`eme 2 est la meilleure possible. Soit $K\in \mathscr{C}^{\infty }(E_{*})$ une fonction strictement positive. On suppose qu'il existe un r\'eel $a\in ]0,1[$ tel que 
$$K(\xi )\leq \frac{a(m-1)}{\Vert \xi \Vert },\ \xi \in E_{*}.$$
Alors il n'existe pas de solution de classe $\mathscr{C}^{2}(\Sigma )$ de l'\'equation $(5.2)$ ci-dessus. En effet, si une telle solution existe, on aura
\begin{equation}
\sum _{n+1\leq \alpha ,\beta \leq n+m-1}B^{\alpha \beta }(u){\tilde
D}_{\alpha \beta }u\geq (m-1)(1+v_{1})-a(m-1)(1+v_{1})^{3/2}.
\end{equation}
En un point $\xi \in \Sigma $ o\`u $u$ atteint son maximum, on a: $\displaystyle v_{1}(\xi )=0$ et dans un rep\`ere $G$-orthonorm\'e diagonalisant $(D_{\alpha \beta }u(\xi ))$, $(5.32)$ s'\'ecrit au point $\xi $ sous la forme 
\begin{equation}
\sum _{n+1\leq \alpha \leq n+m-1}D_{\alpha \alpha }u\geq (m-1)(1-a).
\end{equation}
Or, pour tout $\alpha $, $D_{\alpha \alpha }u(\xi )\leq 0$. Reportons dans $(5.33)$, on voit que $\displaystyle 1\leq a$ ce qui est contradictoire.

\vskip3mm

On obtient la m\^eme conclusion s'il existe un r\'eel $b>1$ tel que
$$K(X)\geq \frac{b(m-1)}{\Vert \xi \Vert },\ \xi \in E_{*}.$$

\vskip4mm

\noindent\textbf{2-} Enfin, on pr\'esente ici un exemple de non-unicit\'e, m\^eme \`a homoth\'etie pr\`es, l'hypoth\`ese de monotonicit\'e $(1.2)$ \'etant satisfaite. Pour cela, consid\'erons le cas o\`u la fonction prescrite est
$$K(\xi )=\frac{m-1}{\Vert \xi \Vert },\ \xi \in E_{*}.$$
Celle-ci v\'erifie l'hypoth\`ese de monotonicit\'e $(1.2)$. En effet, partout dans $E_{*}$, on a :
$$\frac{\partial
\left[\rho K(\rho \xi )\right]}{\partial \rho }=0,\ \mbox{quel que soit }\xi \in \Sigma .$$
La courbure moyenne verticale du fibr\'e $\Sigma _{r}$ est donn\'ee par $K$. D'autre part, si $w\in C^{\infty }(M)$, les calculs men\'es pour prouver le th\'eor\`eme 1 montrent que la courbure moyenne verticale du graphe $Y_{{\tilde w}}$, o\`u ${\tilde w}=w\circ \pi$ est le rel\`evement vertical de $w$ \`a $E_{*}$, est aussi donn\'ee par $K$ et pourtant les hypersurfaces $\Sigma _{r}$ et $Y_{{\tilde w}}$ ne sont pas homoth\'etiques si $w$ n'est pas une constante.

\end{document}